\documentclass[
final
]{dmtcs-episciences}

\usepackage[utf8]{inputenc}
\usepackage{subfigure}

\usepackage[round]{natbib}

\usepackage{amssymb}
\usepackage{amsmath}
\usepackage{amsthm}
\usepackage{algorithm}
\usepackage[noend]{algpseudocode}
\algrenewcommand\algorithmicrequire{\textbf{Input:}}
\algrenewcommand\algorithmicensure{\textbf{Output:}}

\newtheorem{theorem}{Theorem}[section]
\newtheorem{lemma}[theorem]{Lemma} 
\newtheorem{corollary}[theorem]{Corollary}
\newtheorem{proposition}[theorem]{Proposition}

\newtheorem{observation}[theorem]{Observation}

\newcommand{\D}{\ensuremath{\tilde D}}
\newcommand{\lca}{\ensuremath{\mathrm{lca}}}

\author[Huber, Moulton, Scholz]{Katharina T. Huber\affiliationmark{1}
  \and Vincent Moulton\affiliationmark{1}
  \and Guillaume E. Scholz\affiliationmark{2}\thanks{Corresponding author.}}
\title[Arboreal networks and ultrametrics]{Arboreal networks and ultrametrics}

\affiliation{
  School of Computing Sciences, University of East Anglia, Norwich, UK\\
  Institute of Mathematics and Computer Science, University of Greifswald, Greifswald, Germany.
}
\begin{document}

\publicationdata{vol. 28:3}{2026}{25}{10.46298/dmtcs.17508}{2026-02-13; 2026-02-13; 2026-07-29; 2026-10-30}{2026-11-04}

\maketitle
\begin{abstract}
Ultrametrics are an important class of distances used in applications such as 
phylogenetics, clustering and classification theory. Ultrametrics 
are distances that can be represented by an edge-weighted 
rooted tree so that all of the distances in the tree from the root to any leaf 
of the tree are equal.
In this paper, we introduce a generalization of ultrametrics called 
{\em arboreal ultrametrics} which have applications in phylogenetics
and also arise in the theory of distance-hereditary graphs. 
These are partial distances that can be 
represented by an {\em ultrametric arboreal network}, that is, an edge-weighted rooted network whose underlying graph is a tree.
As with ultrametrics, all of the distances in an ultrametric arboreal network from any root to any leaf
below it are equal but, in contrast, the network may have more than one root.
In our two main results we characterize when a partial distance is an 
arboreal ultrametric as well
as proving that, somewhat surprisingly, given any unrooted edge-weighted phylogenetic tree
there is a necessarily 
unique way to insert roots into this tree so as to obtain
an arboreal ultrametric.
\end{abstract}

\section{Introduction}

An {\em arboreal network on $X$} is a digraph whose underlying graph is an undirected or {\em unrooted} phylogenetic tree, and 
whose leaf set is $X$.
Note that such a network can have more than one root, and that if the network
has precisely one root then it is commonly known as a {\em rooted (phylogenetic) tree} \cite{SS03}.
Arboreal networks have applications in evolutionary biology, where the 
leaf set $X$ usually corresponds to a set of species 
and the network represents historical evolutionary relationships between 
the species \cite{HMS22}. Indeed,
there has been growing interest in using networks with multiple roots to represent
evolutionary histories since they can, for example, be used to model evolutionary processes
such as introgression \cite{scholz2019osf} or gene fusion \cite{coleman2015, HJHFLOPWBM13}
where the assumption of a single root can be misleading.
%

\begin{figure}[h]
	\begin{center}
    \includegraphics[scale=0.8]{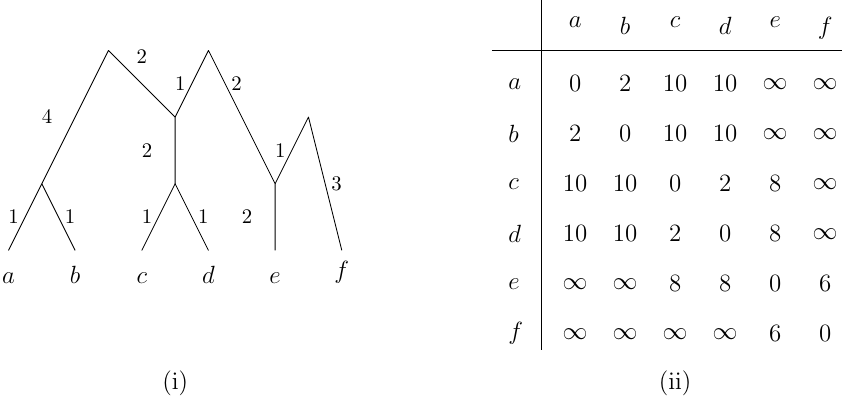}
	\end{center}
	\caption{(i) An ultrametric arboreal network with leaf set $X=\{a,b,\dots,f\}$ 
    and (ii) its associated partial distance, where all
    arcs are directed downwards towards the leaves in $X$.
    For example, the distance from $a$ to $d$ is 10, which is length of the shortest path from $a$ to $d$ that goes via 
    the root that lies above them both, whereas
    the distance from $a$ to $f$ is $\infty$ since there is no root that lies above $a$ and $f$.}
	\label{fig:introduction}
\end{figure}

Despite this interest, the mathematical theory surrounding arboreal networks is still in its infancy. This holds in particular with regards to our understanding of arc-weighted arboreal networks, that is, arboreal networks $N$ in which every arc in $N$ is assigned a non-negative real number called the {\em weight} of that arc. An example of such a network is depicted in Figure~\ref{fig:introduction}(i) where $X=\{a,\ldots, f\}$. By  taking as distance $\tilde{D}(x,y)$ between any two distinct leaves $x,y\in X$ to be the length of the shortest path joining $x$ and $y$ in the underlying tree of an arboreal network $N$ on $X$ in case they are below the same root of $N$ and setting that distance to $\infty$ in case they are not.

 Clearly an arboreal network  induces a special kind of distance on $X$ called a {\em partial distance} -- see Figure~\ref{fig:introduction}(ii) for an example. Formally, such a distance is a symmetric map $\D: X \times X \to \mathbb R_{>0} \cup \{\infty\}$.
 In phylogenetics, partial distances arise when the distance between two taxa is considered irrelevant, or non-existent. 
In the former case (irrelevance), one can consider, given a distance $D$ on $X$ and a threshold $\tau>0$, the partial distance $D_\tau$ on $X$ defined, for all $x,y \in X$ distinct, as $D_\tau(x,y)=D(x,y)$ if $D(x,y) \leq \tau$, and $D_\tau(x,y)=\infty$ otherwise. In the latter case (non-existence), one can assign distance $\infty$ between two taxa that are identified as being unrelated, for example in case of overlapping sets of taxa. We emphasize that, in this context, the equality $\D(x,y)=\infty$ does not necessarily indicate that the distance between $x$ and $y$ is unknown. Indeed, $\D(x,y)=\infty$ is an information about the lack of relation between $x$ and $y$, and thus, shall not be treated as missing information.
This is in contrast to distances which have genuinely missing
entries which have been studied in e.g. \cite{HK16,kettleborough2015reconstructing}.

In some cases the partial distances between the elements on $X$  induced by an arc-weighted arboreal network $N$ on $X$ are such that for any root $\rho$ in $N$ all paths 
from $\rho$ to any leaf in $X$ have the same length, where 
the length of a path is simply the sum of the weights of the arcs in the path (see e.g. Figure~\ref{fig:introduction}(i)). Arc-weighted arboreal networks that enjoy this property are in the center of this paper. To emphasise this, we call such a network an {\em ultrametric arboreal network (on $X$)} and remark in passing that in case the network has a single root, $N$ is also called an {\em equidistant weighted} or {\em ultrametric} tree. Ultrametric
{\em trees} are commonly used in evolutionary biology to represent the evolution of taxa where the length of any 
path from the root to a leaf is assumed to be proportional to the
time that has passed for the root species to evolve to the leaf species \cite[Chapter 9]{Felsenstein2004}.

%
We define an {\em arboreal ultrametric on $X$} to be a partial distance on $X$ if it
can be obtained as the partial distance underlying some ultrametric arboreal network.
Note that 
the arboreal ultrametric $\D$ induced by an ultrametric arboreal network $N$ is an ultrametric if and only if $N$ has a single root.
In addition, as we shall see in Section~\ref{sec:pd}, arboreal ultrametrics are special examples of 
\emph{symbolic arboreal maps} introduced in \cite{huber2024shared}. These maps have close connections with the well-known classes of distance-hereditary graphs \cite{S25} and Ptolemaic graphs \cite{huber2024shared}. As we shall see in Theorem~\ref{thm:diss}, the latter are precisely
the graphs that arise from considering pairs 
of leaves in arboreal networks that have
a common ancestor.

We now give an overview of the main results in this paper. 
One strategy to build ultrametric phylogenetic trees is as follows.
First, construct an unrooted, edge-weighted phylogenetic tree and then, after this, try to insert a root into the tree
so as to create an ultrametric tree (or one close to being ultrametric) using, for example, methods such as the {\em mid-point method} 
or the {\em Farris transform} (see e.g. \cite{dress2007some,kinene2016rooting}). One issue with this strategy is that it 
is not always possible to find such a root. In our first main result we shall show that, in contrast, given 
any phylogenetic tree it is possible to insert roots into the tree so as to create an
ultrametric arboreal network. Moreover, we show that the choice of where to insert the roots is necessarily unique (see Proposition~\ref{pr:uniqueness} and Theorem~\ref{th:rooting}).

We then turn our attention to characterizing when a partial distance is an arboreal ultrametric.
It is well-known (see e.g. \cite[Chapter 7]{SS03}) that a distance $D$ on a set $X$ is an ultrametric if and only if
it satisfies the following strengthening of the metric triangle equality 
$$
D(x,y) \le \max\{D(x,z),D(y,z)\} \mbox{ for all distinct } x,y,z \in X.
$$
In our second main result, we shall give 
a similar characterization for characterizing arboreal ultrametrics.
To decide whether or not a partial distance $\D$ on $X$ is an arboreal ultrametric, it is important
to handle the pairs in $X$ for which $\D$ is infinity.
To do this, we consider the graph $G_{\D}$ that has vertex set $X$ in which two distinct vertices $x$ and $y$ are joined by an edge if $\D(x,y)<\infty$. We then prove in Theorem~\ref{th:ultra} that a partial distance $\D$ is an
arboreal ultrametric if and only if (i) $G_{\D}$  is a connected, chordal graph, (ii) $\D$ satisfies 
the above 3-point ultrametric condition in case $\D$ is defined for all pairs in the triple, and (iii) 
$\D$ satisfies an additional 4-point condition. The 
4-point condition is given in full in Theorem~\ref{th:ultra}.

The rest of this paper is organized as follows.
In Section~\ref{sec:prlm}, we present some preliminaries. Then in Section~\ref{sec:au}
we prove the aforementioned result about rooting trees 
to obtain ultrametric arboreal networks (Theorem~\ref{th:rooting}). In Section~\ref{sec:pd}
we present our characterization for 
arboreal ultrametrics (Theorem~\ref{th:ultra}), before concluding with
a brief discussion of some future directions in Section~\ref{sec:future}.

\section{Preliminaries}\label{sec:prlm}

We shall assume throughout the paper that $X$ is a finite set for which $|X|\geq 2$ holds.

\subsection*{Graphs}

A graph $G$ is an ordered pair $(V,E)$, where $V=V(G)$ is a finite set of elements, called \emph{vertices} (of $G$), and $E=E(G)$ is a set of pairs of distinct elements of $V$. If the pairs in $E$ are not ordered, we call them \emph{edges}, and we say that $G$ is \emph{undirected}. We denote an edge between two vertices $u$ and $v$ by $\{u,v\}$.  If the pairs in $E$ are ordered, we call them \emph{arcs}, and we say that $G$ is \emph{directed}. For two vertices $u$ and $v$ of $V$ we denote the arc from $u$ to $v$ by $(u,v)$. For an arc $(u,v)$ of a directed graph $G$, we say that $u$ is a \emph{parent} of $v$, and $v$ is a \emph{child} of $u$.

A \emph{path} in an undirected (\emph{resp.} directed) graph $G$ is a sequence $x_1, \ldots, x_k$, $k \geq 1$ of pairwise distinct elements of $X$ such that for all $i \in \{1,\ldots,  k-1\}$, $\{x_i,x_{i+1}\}$ is an edge of $G$ (\emph{resp.} $(x_i,x_{i+1})$ is an arc of $G$). The \emph{length} of a path is the number of its edges (\emph{resp.} arcs). More specifically, a path $x_1, \ldots, x_k$, $k \geq 1$ has length $k-1$. A \emph{cycle} is a sequence $x_1, \ldots, x_k$, $k \geq 3$ of elements of $G$ such that $x_1, \ldots, x_k$ is a path of $G$, and in addition, $\{x_k,x_1\}$ is an edge of $G$ (\emph{resp.} $(x_k,x_1)$ is an arc of $G$). As in the case of paths, the \emph{size} of a cycle is the number of edges (\emph{resp.} arcs) it contains. A graph that does not contain any cycle is called \emph{acyclic}. In the undirected case, we also sometimes refer to connected acyclic graphs as \emph{trees}.

We say that two undirected graphs $G=(V,E)$ and $G'=(V',E')$ are \emph{isomorphic} if there exists a bijection $\phi: V \to V'$ such that $\{u,v\} \in E$ if and only if $\{\phi(u),\phi(v)\} \in E'$. This definition naturally extends to directed graphs.

\begin{figure}[h]
	\begin{center}
		\includegraphics[scale=0.8]{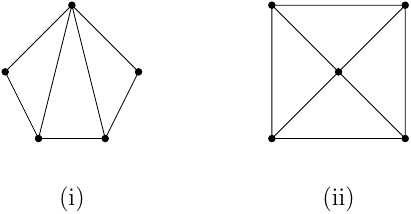}
	\end{center}
	\caption{(i) The gem is an example of a chordal graph on $5$ vertices, and it is the only chordal forbidden induced subgraphs for Ptolemaic graphs. (ii) The wheel $W_5$ is a non-chordal graph that is obtained from the gem by  adding a single edge.}
	\label{fig:gw}
\end{figure}

Given an undirected graph $G=(V,E)$ and a non-empty subset $Y \subseteq V$, the \emph{subgraph of $G$ induced by $Y$}, denoted by $G[Y]$, is the graph with vertex set $Y$, and with edge set the set $\{\{u,v\} \in E\,:\, u,v \in Y\}$. We say that an undirected graph $G$ is \emph{chordal} if it does not contain a cycle of size $4$ or more as an induced subgraph. Among the undirected graphs of interest to us are the {\em gem} which is a path $P$ of length $3$ together with a further vertex $x$ not on $P$ that is adjacent to all vertices of $P$ (see Figure~\ref{fig:gw}(i)), and the \emph{wheel} $W_k$, $k \geq 4$, which is a cycle $C$ of length $k-1$ together with a vertex $x$ adjacent to all vertices of $C$ (see Figure~\ref{fig:gw}(ii) for the wheel $W_5$) . We say that a graph $G$ is \emph{Ptolemaic} \cite{H81} if $G$ is chordal and does not contain the gem as an induced subgraph.

Three key operations on graphs will be of interest to us throughout this contribution. Let $G=(V,E)$ be an undirected graph, $v$ a vertex of $G$ of degree $2$, and let $u$ and $w$ be the vertices adjacent to $v$ in $G$. Then, the first operation is the \emph{suppressing} of $v$ operation which consists of removing $v$ and the edges $\{u,v\}$ and $\{v,w\}$ from $G$, and adding the edge $\{u,w\}$. If $G$ is directed, and $v$ is a vertex of $G$ with a unique parent $u$ and a unique child $w$, then \emph{suppressing} of $v$ is the operation that consists of removing $v$ and the arcs $(u,v)$ and $(v,w)$ from $G$, and adding the arc $(u,w)$.

The second operation is an operation that reverses the suppressing operation. More precisely,  given an edge $\{u,w\}$ of an undirected graph $G$, then \emph{subdividing} of $\{u,w\}$ consists of removing the edge $\{u,w\}$ from $G$, adding a new vertex $v$, and adding the edges $\{u,v\}$ and $\{v,w\}$. If $G$ is directed and $(u,w)$ is an arc of $G$, then the \emph{subdividing} of $(u,w)$ operation consists of removing the arc $(u,w)$ from $G$, adding a new vertex $v$, and adding the arcs $(u,v)$ and $(v,w)$. In either case, we call such a vertex $v$ a {\em subdivision vertex}.

Our final operation is the \emph{contracting} of an edge $\{u,v\}$ operation which in an undirected graph $G$ consists of removing $v$ from $G$, and replacing all edges $\{v,w\}$, $w \neq u$ of $G$ with the edge $\{u,w\}$. If $G$ is a directed graph, then the  \emph{contracting} of an arc $(u,v)$ operation consists of removing $v$ from $G$, replacing all arcs $(v,w)$ of $G$ with the arc $(u,w)$, and replacing all arcs $(w,v)$, $w \neq u$ of $G$ with the arc $(w,u)$.

\subsection*{Networks}

A \emph{network $N$ (on $X$)} 
is a connected directed acyclic graph with leaf set $X$  
such that all vertices of indegree 0 have 
outdegree at least 2, all vertices of outdegree 0 have indegree 1, and no 
vertices have indegree and outdegree equal to 1. We call the vertices of indegree 0 the {\em roots} of a network and the vertices with outdegree $0$ its {\em leaves}. We denote by $L(N)$ the set of leaves of $N$ and by $R(N)$ the set of roots of $N$. Also, we
	put $r(N)=|R(N)|$. Note that we must have $r(N)\geq 1$. A network with a single root is commonly called  a {\em rooted phylogenetic network}. In this case, if the
    undirected graph obtained by ignoring directions in the network is
    a tree, it is called a {\em rooted phylogenetic tree} (see e.g. \cite{S16} for more details on such networks).
We say that two networks $N$, $N'$ on $X$ are \emph{isomorphic} if there exists a digraph isomorphism $\phi$ from the vertex set $V(N)$ of $N$ to the vertex set $V(N')$ of $N'$ that is the identity on $X$.

Assume for the remainder of this section that  $N$ is a network on $X$.  For $u,v$ two vertices of $N$, we say that $u$ is an \emph{ancestor} of $v$ if there is a path from $u$ to $v$ in $N$. In this case, we also refer to $v$ as a \emph{descendant} of $u$. If in addition, $u \neq v$, then we say that $u$ is a \emph{proper ancestor} of $v$ and that $v$ a \emph{proper descendant of $u$}. 

\begin{figure}[h]
	\begin{center}
		\includegraphics[scale=0.66]{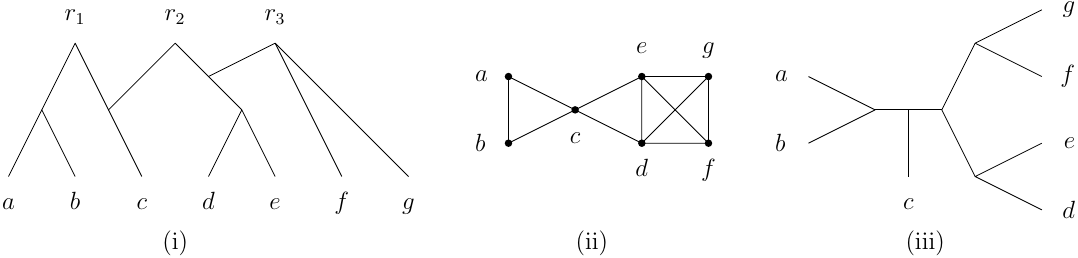}
	\end{center}
	\caption{(i) An arboreal network $N$ on $X=\{a,b,c,d,e,f,g\}$, with three roots $r_1, r_2$ and $r_3$. (ii) The shared ancestry graph $\mathcal A(N)$ of $N$. (iii) The underlying graph $\overline N$ of $N$.}
	\label{fig:prlm}
\end{figure}

For $v$ a vertex of $N$, we denote by $C_N(v)$ the set of all leaves of $N$ that are a descendant of $v$. We say that two distinct leaves $x,y \in X$ \emph{share an ancestor} in $N$ if there exists a vertex $v$ of $N$ such that $x,y \in C_N(v)$. This notion allows us to define the \emph{shared-ancestry graph} $\mathcal A(N)$ of $N$ as follows. The vertex set of $\mathcal A(N)$ is $X$, and two distinct elements $x,y \in X$ are joined by an edge if and only if $x$ and $y$ share an ancestor in $N$. As an example, for $N$ the network on $X=\{a,b,c,d,e,f,g\}$ depicted in Figure~\ref{fig:prlm}(i), we present the shared-ancestry graph $\mathcal A(N)$ of $N$ in Figure~\ref{fig:prlm}(ii).

If two distinct leaves $x, y \in X$ share an ancestor in $N$, we say that a vertex $v$ is a \emph{lowest common ancestor} of $x$ and $y$ if $v$ is an ancestor of both $x$ and $y$, and no child of $v$ is also an ancestor of both $x$ and $y$. By definition, two vertices sharing an ancestor in $N$ have at least one lowest common ancestor. However, that ancestor may not be unique.
If two leaves $x$ and $y$ of $N$ do have a unique lowest common ancestor,
	then we denote it by $\lca_N(x,y)$.

Let $R_2(N)$ be the set of roots of $N$ with outdegree  $2$. We 
define the \emph{underlying graph $\overline N$} of $N$ as the undirected graph with vertex set $V(\overline N)=V(N)\setminus R_2(N)$, and with edge set those pairs $\{u,v\}$ of vertices in $V(\overline N)$ such that one of the following holds:
\begin{itemize}
\item[-] One of $(u,v)$ or $(v,u)$ is an arc of $N$.
\item[-] There exists a root $r \in R_2(N)$ such that $(r,u)$ and $(r,v)$ are arcs of $N$.
\end{itemize}

In other words,  $\overline N$ is the undirected graph obtained from $N$ by suppressing the direction of the arcs and suppressing resulting vertices of degree $2$. In particular, the leaf set of $\overline N$ is $X$. As an example, for $N$ the network depicted in Figure~\ref{fig:prlm}(i), the underlying undirected graph $\overline N$ of $N$ is depicted in Figure~\ref{fig:prlm}(iii).

Following \cite{HMS22}, we say that a network $N$ on $X$ is \emph{arboreal} if its underlying graph $\overline N$ is a {\em phylogenetic tree (on $X$)}, that is, an unrooted tree which does not have any vertices of degree two and whose leaf set is $X$. Note that the network $N$ depicted in Figure~\ref{fig:prlm}(i) is arboreal, since the underlying graph $\overline N$ of $N$, depicted in Figure~\ref{fig:prlm}(iii), is a phylogenetic tree on $L(N)$. Note that a network $N$ is arboreal if and only if for all arcs $a$ of $N$, the removal of $a$ from $N$ disconnects $N$. Also, note that the least common ancestor of a pair of leaves in an arboreal network, if it exists, is 
necessarily unique \cite[Proposition 7.1]{huber2024shared}, and that the shared-ancestry graph of an arboreal network 
is always Ptolemaic \cite[Proposition 6.3]{huber2024shared}.

For $T$ a phylogenetic tree on $X$, an arboreal network $N$ such that $\overline N$ and $T$ are isomorphic can be obtained by subdividing some of the edges of $T$, and assigning directions to the edges of the resulting graph in such a way that all subdivision vertices have indegree $0$, and a vertex $v$ has outdegree $0$ if and only if it is a leaf of $T$. Note that in view of the first requirement, an edge of $T$ cannot be subdivided more than once in this process. Note also that a non-leaf vertex that is not a subdivision vertex is allowed to have indegree $0$. We call a network obtained this way an \emph{uprooting} of $T$. It is straight-forward to check that, although 
the underlying tree of an arboreal network is uniquely defined, there may in general be several non-isomorphic uprootings of a given phylogenetic tree.

%
%

\subsection*{Distances and weighted networks}

Let $G$ be an undirected graph. A \emph{weighting} of $G$ is a map $\lambda: E(G) \to \mathbb R_{>0}$ that assigns to each edge of $G$ a positive value. We call the pair $(G,\lambda)$ a \emph{weighted graph}. For $e \in E(G)$, we call $\lambda(e)$ the \emph{length} of $e$. We extend this definition to directed graphs in the obvious way,
and will tend to use the symbol $\omega$ instead of $\lambda$ to designate weightings in such graphs.
In addition, the notion of isomorphism can be generalized to weighted graphs as follows. We say that two weighted graphs $(G,\lambda)$ and $(G',\lambda')$ are \emph{isomorphic} if $G$ and $G'$ are isomorphic
via a map $\phi:V(G) \to V(G')$, and the bijection $\phi$ satisfies $\lambda(\{u,v\})=\lambda'(\{\phi(u),\phi(v)\})$ for all edges $\{u,v\}$ of $G$. Again, this definition extends naturally to directed graphs.

Let $(T,\lambda)$ be a weighted phylogenetic tree on $X$. Since $T$ is a tree, then for any pair $u,v$ of vertices of $T$, there exists a unique path $P_T(u,v)$ between $u$ and $v$ in $T$. The \emph{distance} $l_{(T,\lambda)}(u,v)$ between $u$ and $v$ in $(T,\lambda)$ is defined as the sum of the lengths of all edges lying on $P_T(u,v)$. In particular, $l_{(T,\lambda)}(u,v) \geq 0$ always holds, with equality holding if and only if $u=v$. Viewing $l_{(T,\lambda)}$ as a map from $V(T) \times V(T)$ into the non-negative reals, then the restriction of the map $l_{(T,\lambda)}$ to $X \times X$ induces a map $D_{(T,\lambda)}: X \times X \to \mathbb R_{\geq 0}$ defined by putting $D_{(T,\lambda)}(x,y)=l_{(T,\lambda)}(x,y)$, for all $x,y\in X$.

Now, consider an arbitrary distance $D$ on $X$, that is, a map $D: X \times X \to \mathbb R_{\geq 0}$, that is, 
a map  such that $D$ is \emph{symmetric} (i.e. $D(x,y)=D(y,x)$ for all $x,y \in X$), and 
$D$ vanishes precisely on the diagonal (i.e. $D(x,y)=0$ if and only if $x=y$).
We say that $D$ is \emph{tree-like} if there exists a weighted phylogenetic tree $(T,\lambda)$ on $X$ such that $D=D_{(T,\lambda)}$. In that case, we say that $(T,\lambda)$ \emph{represents} $D$. As is well-known (\cite{B71}, see also \cite[Section 7.1]{SS03}), we have:

\begin{theorem}\label{th:treelike}
Let $D$ be a distance on $X$. Then $D$ is tree-like if and only if  $D$ satisfies the {\em four-point condition}, that is, for all (not necessarily distinct) $x,y,z,u \in X$, $D(x,y)+D(z,u) \leq \mathrm{max}\{D(x,z)+D(y,u),D(x,u)+D(y,z)\}$ holds. Moreover, if $D$ is tree-like, then there exists a unique (up to isomorphism) phylogenetic tree $T$ and a unique weighting $\lambda$ of $T$ such that $D=D_{(T,\lambda)}$.
\end{theorem}

Consider now a weighted arboreal network $(N,\omega)$. The \emph{weighted underlying phylogenetic tree} of $(N,\omega)$ is the weighted phylogenetic tree $(\overline N,\overline{\omega})$, where $\overline{\omega}$ is defined, for all edges $\{u,v\}$ of $\overline N$, as:
\begin{itemize}
\item[-] $\overline{\omega}(\{u,v\})=\omega((u,v))$ (\emph{resp.} $\omega((v,u))$) if $(u,v)$ (\emph{resp.} $(v,u)$) is an arc of $N$.
\item[-] $\overline{\omega}(\{u,v\})=\omega((r,u))+\omega((r,v))$ if $u$ and $v$ are children of some $r \in R_2(N)$ in $N$.
\end{itemize}

\noindent Note that $(\overline N,\overline{\omega})$ is uniquely determined by $(N,\omega)$. In view of this, and since $\overline N$ is a phylogenetic tree, we define $D_{(N,\omega)}:=D_{(\overline N,\overline{\omega})}$, and we call $D_{(N,\omega)}$ the \emph{distance induced by} $(N,\omega)$.

We now define the concept of an ultrametric arboreal network, a concept that generalises the 
definition of an ultrametric tree. To do this, we require further notation.
Let $(N,\omega)$ be a weighted arboreal network on $X$. For $u,v$ two vertices of $N$ such that $u$ is an ancestor of $v$ in $N$, there is a unique directed path $P_N(u,v)$ from $u$ to $v$ in $N$ (see e.g. \cite[Lemma~3]{S25}). We denote the sum of the lengths of all arcs of $P_N(u,v)$  by $l_{(N,\omega)}(u,v)$. Note that if $u \notin R_2(N)$, then by definition $u$ and $v$ are vertices of $\overline N$, so $l_{(N,\omega)}(u,v)=l_{(\overline N,\overline{\omega})}(u,v)$ holds.
We say that a weighted arboreal network $(N,\omega)$ is \emph{ultrametric} if for all vertices $u$ of $N$, and all leaves $x,y \in C_N(u)$, we have $l_{(N,\omega)}(u,x)=l_{(N,\omega)}(u,y)$.
This definition generalises the notion of an \emph{ultrametric tree} 
(sometimes also known as an \emph{equidistant tree}, see e.g. \cite[Section 7.2]{SS03}) which, in our terminology, 
is simply an ultrametric arboreal network $N$ with a single root.

We conclude this section by recalling the concept of an ultrametric.
We say that a distance $D: X \times X \to \mathbb R_{\geq 0}$ is an \emph{ultrametric}
if $D = D_{(T,\omega)}$ for $(T,\omega)$ an ultrametric tree.
In this case, we say that $(T,\omega)$ \emph{represents} $D$. Similar to the case of 
tree-like distances, ultrametrics can be characterized in terms of the following 3-point condition 
(see e.g. \cite[Theorem 7.2.5]{SS03}):

\begin{theorem}\label{th:ultree}
Let $D$ be a distance on $X$. Then $D$ is an ultrametric if and only if $|X| \leq 2$ or, for all $x,y,z \in X$ distinct, $D(x,y) \leq \mathrm{max}\{D(x,z),D(y,z)\}$. Moreover, if this holds, then there exists a unique (up to isomorphism) ultrametric tree $(T,\omega)$ such that $D=D_{(T,\omega)}$.
\end{theorem}

\section{Ultrametric uprootings}\label{sec:au}

\begin{figure}[h]
	\begin{center}
    \includegraphics[scale=0.8]{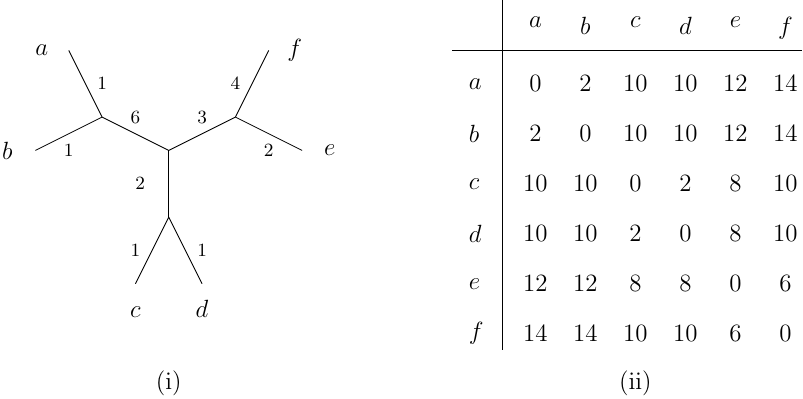}
	\end{center}
	\caption{(i) A weighted phylogenetic tree $(T,\lambda)$ that admits the ultrametric network in Figure~\ref{fig:introduction} as an ultrametric uprooting. (ii) The distance $D_{(T,\lambda)}$.}
	\label{fig:uprooting}
\end{figure}

Often in phylogenetics, biologists first compute a weighted phylogenetic tree $(T,\lambda)$
for their data and then insert a root into this tree 
to obtain a rooted phylogenetic tree in which the root represents the last
common ancestor of all of the leaves.
For most data sets it is not possible to do this in such a way 
that the resulting rooted tree is an ultrametric tree, since this would be
equivalent to $D_{(T,\lambda)}$ being an an ultrametric, which is often not the case. 
For example, for $(T,\lambda)$ the weighted phylogenetic tree depicted in Figure~\ref{fig:uprooting}(i), the distance $D:=D_{(T,\lambda)}$ (Figure~\ref{fig:uprooting}(ii)) is not an ultrametric, 
since we have $D(a,e)=12>10=\mathrm{max}\{D(a,c),D(c,e)\}$. 

Even so, if one were trying to build an arboreal network instead of a rooted tree, 
it could still be of interest to understand whether more
roots might be inserted into a weighted phylogenetic tree
order to obtain an ultrametric arboreal network
even if no ultrametric rooting of the tree were possible.
In this section, we shall show that it is in fact {\em always} possible to insert some roots into 
a given weighted phylogenetic tree to obtain such a network and, in fact, for any 
given weighted phylogenetic tree there is only one way to do this (Theorem~\ref{th:rooting}).

\subsection*{Definitions}

For $(T,\lambda)$ a weighted phylogenetic tree with leaf set $X$, we 
say that a weighted arboreal network $(N,\omega)$ on $X$ is a \emph{weight-preserving uprooting} of $(T,\lambda)$ if 
$(T,\lambda)$ is isomorphic to the underlying weighted graph $(\overline N,\overline{\omega})$ of $N$.
In addition, we say that $(N,\omega)$ is an \emph{ultrametric uprooting} of $(T,\lambda)$ 
if $(N,\omega)$ is a weight-preserving uprooting of $(T,\lambda)$, which, in addition, is ultrametric. 
For example, the ultrametric network depicted in Figure~\ref{fig:introduction}(i) is 
an ultrametric uprooting of the weighted phylogenetic tree depicted in Figure~\ref{fig:uprooting}(i).

Note that, as we have seen above, there exists in general more than one uprooting of a phylogenetic tree $T$. 
Similarly, for a given uprooting $N$ of $T$, there are 
infinitely many edge-weightings $\omega$ of $N$ such that $\lambda=\overline{\omega}$. 
To see this, let $(N,\omega)$ be a weight-preserving uprooting of $(T,\lambda)$. 
By definition, $\omega((u,v))=\lambda(\{u,v\})$ for all arcs $(u,v)$ of $N$ 
for which $u \notin R_2(N)$. However, if $u \in R_2(N)$, the only 
requirement on $\omega((u,v))$ is that the equality $\omega((u,v))+\omega((u,v'))=\lambda(\{v,v'\})$ 
holds for $v'$ the second child of $u$ in $v$.

We start with proving 
that any weighted phylogenetic tree admits an ultrametric uprooting (Proposition~\ref{pr:uur}). Since our proof is constructive we shall present an algorithm called \textsc{Ultrametric Uprooting} in Algorithm~\ref{alg:ulup}
    and prove that it always produces such an uprooting (Proposition~\ref{pr:uur}). We first introduce some additional notation concerning phylogenetic trees.

Let $T$ be a phylogenetic tree on $X$. 
Given a subset $Y \subseteq X$, we denote by $T|_Y$ the phylogenetic
tree obtained from $T$ by removing all leaves in $X \setminus Y$ as well as 
all new leaves created in the process, and suppressing all resulting vertices of degree $2$.
Note that all edges $e=\{u,v\}$ in $T|_Y$ coincide with the path $P_T(u,v)$ between $u$ and $v$ in $T$. 
In view of this, if $T$ is equipped with a weighting $\lambda$, we define 
the weight $\lambda|_Y(e)$ of an edge $e$ of $T|_Y$ as the sum of the 
weights (under $\lambda$) of the edges in $P_T(u,v)$.
In addition, we say that a leaf $x$ of $T$ is \emph{in a cherry} if 
the (necessarily unique) vertex $v$ adjacent to $x$ in $T$ is adjacent 
to a leaf $y \in X$ that is distinct from $x$. 
If in addition, $T$ is equipped with a weighting $\lambda$, and $x \in X$ 
is such that $\lambda(\{v,x\}) \geq \lambda(\{v,y\})$ for 
all leaves $y \in X$ adjacent to $v$ in $T$, then we say that $x$ 
is a \emph{long-end of a cherry} in $T$. 

\begin{figure}[h]
\begin{center}
\includegraphics[scale=0.6]{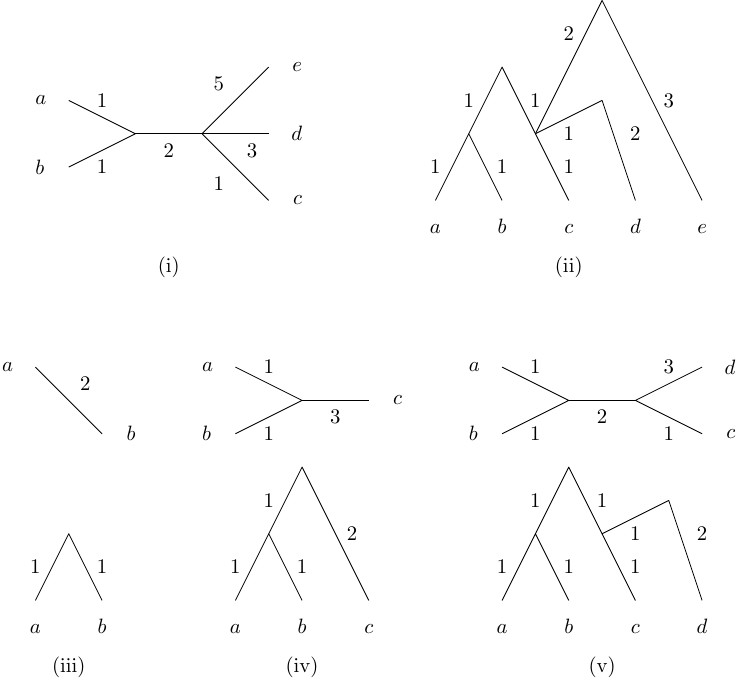}
\caption{(i) A weighted phylogenetic tree $(T,\lambda)$ on $X=\{a,b,c,d,e\}$. (ii) the output of Algorithm~\ref{alg:ulup} applied on $(T,\lambda)$ and the wcps $a,b,c,d,e$. (iii), (iv) and (v) Top, the trees $T|_{\{a,b\}}$, $T|_{\{a,b,c\}}$ and $T|_{\{a,b,c,d\}}$, respectively. Bottom, the intermediate steps of the algorithm. See text for details.}
\label{fig:ex1}
\end{center}
\end{figure}

We say that an ordering $x_1, \ldots, x_n$, $n=|X|$ of the elements of $X$ is 
a \emph{cherry-picking sequence} \cite{HLS13}\footnote{Note that, in \cite{HLS13} a cherry picking sequence is defined on a set $S$ of rooted phylogenetic trees. Here, we restrict to the case $|S|=1$.} (or \emph{cps} for short) if, 
for all $i \in \{3, \ldots, n\}$, $x_i$ is part of a cherry in $T|_{\{x_1, \ldots, x_i\}}$. 
If $T$ is equipped with a weighting $\lambda$, we call a 
cherry-picking sequence $x_1, \ldots, x_n$ a \emph{weighted cherry-picking sequence} 
(or \emph{wcps} for short) if, for all $i \in \{3, \ldots, n\}$, $x_i$ 
is the long end of a cherry in $T|_{\{x_1, \ldots, x_i\}}$. 
As an example, consider the weighted phylogenetic tree $(T,\lambda)$ on $X=\{a,b,c,d,e\}$ 
depicted in Figure~\ref{fig:ex1}(i). Then, the sequences 
$a,b,c,d,e$ and $e,d,c,b,a$ are both cherry picking sequences. 
However, only the former is a weighted cherry picking sequence. 
Indeed, to see that $e,d,c,b,a$ is not a wcps, we remark that, even 
though $c$ is part of a cherry in $T|_{\{e,d,c\}}$, it is 
not the long-end of a cherry, as the (unique) long-end of that cherry is $e$. Note that the same also holds for $b$ in $T|_{\{e,d,c,b\}}$.

As remarked in \cite{D19}, all phylogenetic trees $T$ on $X$ 
admit a cps, which can be computed as follows. Put $n=|X|$.
First, consider a non-leaf vertex $v$ of $T$ adjacent to at least two leaves. 
Note that in a phylogenetic tree, such a vertex always exists. 
Then choose $x_n$ as one of the leaves adjacent to $v$ in $T$, and repeat 
this process with $T$ replaced by $T|_{X \setminus \{x_n\}}$. After $n-2$ instances, 
we are left with a tree that has two leaves $x$ and $x'$, which can be independently defined 
as $x_1$ and $x_2$. Unsurprisingly, this approach can also be used to 
compute a wcps for some weighted phylogenetic tree $(T,\lambda)$. Indeed, 
after choosing a non-leaf vertex $v$ of $T$ adjacent to at least two leaves, one can pick a leaf $x$ adjacent to $v$ 
such that $\lambda(\{v,x\}) \geq \lambda(\{v,x'\})$ for all 
leaves $x'$ adjacent to $v$. By definition, this results in a wcps.

\subsection*{Algorithm}

We are now ready to present our algorithm, \textsc{Ultrametric Uprooting}, for uprooting a weighted phylogenetic tree  (Algorithm~\ref{alg:ulup}).
Before proving that this algorithm is correct, 
we illustrate its inner-workings using a small example.
Consider the weighted phylogenetic tree $(T,\lambda)$ on $X=\{a,b,c,d\}$ depicted in 
Figure~\ref{fig:ex1}. As remarked above, $a,b,c,d,e$ is a wcps 
of $(T,\lambda)$. So assume that $a,b,c,d,e$ is the wcps computed at Line~\ref{l:wcps}.

\begin{algorithm}
  \caption{\textsc{Ultrametric Uprooting}}
  \label{alg:ulup}
  \begin{algorithmic}[1]
     \Require  A weighted phylogenetic tree $(T,\lambda)$ with leaf set $X$
    \Ensure   An ultrametric uprooting $(N,\omega)$ of $(T,\lambda)$
    \State Compute a wcps $x_1, \ldots, x_{|X|}$ \label{l:wcps}
    \State Initialize $N$ as the network with vertex set $\{x_1,x_2,r_2\}$ and arc set $\{(r_2,x_1),(r_2,x_2)\}$ \label{l:N0}
    \State Define $\omega((r_2,x_1))=\omega((r_2,x_2))=\frac{1}{2}\lambda_2(\{x_1,x_2\})$  \label{l:o0}
    \For{$i$ from $3$ to $|X|$}\label{l:fori}
       \State Put $T_i=T|_{\{x_1, \ldots, x_i\}}$ and $\lambda_i=\lambda|_{\{x_1, \ldots, x_i\}}$
       \State Let $v_i$ be the vertex of $T_i$ adjacent to $x_i$, and let $x$ be a further leaf of $T_i$ adjacent to $v_i$ \label{l:defvx}
       \State Let $p_x$ be the parent of $x$ in $N$ \label{l:defvp}
       \If{$\omega((p_x,x)) > \lambda_i(\{v_i,x\})$} \label{l:ifbig}
         \State Subdivide $(p_x,x)$ by introducing a new vertex $v$ \label{l:vbig}
         \State Define $\omega((v,x))=\lambda_i(\{v_i,x\})$ and $\omega((p_x,v))=\omega((p_x,x))-\lambda_i(\{v_i,x\})$ \label{l:ombig}
      \ElsIf{$\omega((p_x,x)) < \lambda_i(\{v_i,x\})$} \Comment{$p_x$ has outdegree $2$ (see proof)} \label{l:ifsm} 
          \State Let $u$ be the child of $p_x$ distinct from $x$
          \If{$\lambda_i(\{v_i,x\})=\omega((p_x,x))+\omega((p_x,u))$}
            \State{Define $v=u$} \label{l:vsm0}
          \Else
          \State Subdivide $(p_x,u)$ by introducing a new vertex $v$ \label{l:vsm}
          \State Define $\omega((p_x,v))=\lambda_i(\{v_i,x\})-\omega((p_x,x))$ and $\omega((v,u))=\omega((p_x,u))+\omega((p_x,x))-\lambda_i(\{v_i,x\})$ \label{l:omsm}
          \EndIf
       \Else \Comment{$\omega((p_x,x)) = \lambda_i(\{v_i,x\})$} \label{l:ifeq}
         \State Define $v=p_x$  \label{l:veq}
       \EndIf
       \If {$\lambda_i(\{v_i,x_i\})=l_{(N,\omega)}(v,z)$ for some $z \in C_N(v)$}\label{l:ifdown}
         \State Add to $N$ the arc $(v,x_i)$ \label{l:vdown}
         \State Define $\omega((v,x_i))=\lambda_i(\{v_i,x_i\})$ \label{l:omdown}
       \Else \label{l:iftop}
          \State Add to $N$ a new vertex $r$ and the arcs $(r,v)$ and $(r,x_i)$ \label{l:vtop}
         \State Define $\omega((r,v))=\frac{1}{2}(\lambda_i(\{v_i,x_i\})-l_{(N,\omega)}(v,z))$ and $\omega((r,x_i))=\frac{1}{2}(\lambda_i(\{v_i,x_i\})+l_{(N,\omega)}(v,z))$ \label{l:omtop}
       \EndIf
    \EndFor  
\State \Return $(N,\omega)$ \label{l:end}
  \end{algorithmic}
\end{algorithm}

First, we initialize $N$ as the network with a single root and two leaves $a$ and $b$ (Line~\ref{l:N0}). Since $\lambda_2(\{a,b\})=2$, we assign length $\lambda_2(\{a,b\})/2=1$ to both arcs of $N$ (Line~\ref{l:o0}). This results in the network depicted in Figure~\ref{fig:ex1}(iii), bottom, which is an ultrametric uprooting of the weighted phylogenetic tree depicted in Figure~\ref{fig:ex1}(iii), top. Next, we enter the loop at Line~\ref{l:fori}.

In the first instance of the loop, we consider the leaf $x_3=c$ of $T|_{\{a,b,c\}}$ (Figure~\ref{fig:ex1}(iv), top). The vertex $v_3$ adjacent to $c$ in $T|_{\{a,b,c\}}$ is adjacent to both $a$ and $b$, so we can independently choose $x=a$ or $x=b$ at Line~\ref{l:defvx}. Suppose we choose $x=a$. We have $\lambda_3(\{v_3,a\})=1$ and $\omega((p_a,a))=1$, so $\omega((p_a,a))=\lambda_3(\{v_i,a\})$, and neither of the conditions at Lines~\ref{l:ifbig} and \ref{l:ifsm} ar satisfied. Hence we enter the subcase at Line~\ref{l:ifeq}, and we put $v=p_a$ (Line~\ref{l:veq}). Next, we have $\lambda_3(\{v_3,c\})=3$ and $l_{(N,\omega)}(v,a)=1$. Since $a \in C_N(v)$, we are in the situation of Line~\ref{l:iftop}. Following Lines~\ref{l:vtop} and \ref{l:omtop}, we therefore add a new vertex $r$ and the arcs $(r,c)$ and $(r,v)$, and we put $\omega((r,v))=\frac{1}{2}(\lambda_3(\{v_3,c\})-l_{(N,\omega)}(v,a))=1$ and $\omega((r,c))=\frac{1}{2}(\lambda_3(\{v_3,c\})+l_{(N,\omega)}(v,a))=2$. This results in the network depicted in Figure~\ref{fig:ex1}(iv), bottom, which is an ultrametric uprooting of the weighted phylogenetic tree depicted in Figure~\ref{fig:ex1}(iv), top.

In the second instance of the loop, we consider the leaf $x_4=d$ of $T|_{\{a,b,c,d\}}$ (Figure~\ref{fig:ex1}(v), top). The vertex $v_4$ adjacent to $d$ in $T|_{\{a,b,c,d\}}$ is adjacent to $c$ only, so we must choose $x=c$ at Line~\ref{l:defvx}. We have $\lambda_4(\{v_4,c\})=1$ and $\omega((p_c,c))=2$, so $\omega((p_c,c))>\lambda_4(\{v_4,c\})$. Hence, the condition of Line~\ref{l:ifbig} holds. We then subdivide the arc $(p_c,c)$ by introducing a new vertex $v$(Line~\ref{l:vbig}), and we put $\omega((v,c))=\lambda_4(\{v_4,c\})=1$ and $\omega((p_c,v))=\omega((p_c,c))-\lambda_4(\{v_4,c\})=1$ (Line~\ref{l:ombig}). Next, we have $\lambda_4(\{v_4,d\})=3$ and $l_{(N,\omega)}(v,c)=1$. Since $c \in C_N(v)$, we are in the situation of Line~\ref{l:iftop}. Following Lines~\ref{l:vtop} and \ref{l:omtop}, we therefore add a new vertex $r$ and the arcs $(r,d)$ and $(r,v)$, and we put $\omega((r,v))=\frac{1}{2}(\lambda_4(\{v_4,d\})-l_{(N,\omega)}(v,c))=1$ and $\omega((r,d))=\frac{1}{2}(\lambda_4(\{v_4,d\})+l_{(N,\omega)}(v,c))=2$. This gives rise to the network depicted in Figure~\ref{fig:ex1}(v), bottom, which is an ultrametric uprooting of the weighted phylogenetic tree depicted in Figure~\ref{fig:ex1}(v), top.

In the third and last instance of the loop, we consider the leaf $x_5=e$ of $T|_{\{a,b,c,d,e\}}=T$ (Figure~\ref{fig:ex1}(i)). The vertex $v_5$ adjacent to $e$ in $T$ is adjacent to both $c$ and $d$, so we 
can independently choose $x=c$ or $x=d$ at Line~\ref{l:defvx}. Suppose we choose $x=c$. We have $\lambda_5(\{v_5,c\})=1$ and $\omega((p_c,c))=1$, so $\omega((p_c,c))=\lambda(\{v_5,c\})$, and neither of the conditions at Lines~\ref{l:ifbig} and \ref{l:ifsm} ar satisfied. Hence, we enter the subcase at Line~\ref{l:ifeq}, and we put $v=p_c$ (Line~\ref{l:veq}). Next, we have $\lambda_5(\{v_5,e\})=5$ and $l_{(N,\omega)}(v,c)=1$. Since $c \in C_N(v)$, we are in the situation of Line~\ref{l:iftop}. Following Lines~\ref{l:vtop} and \ref{l:omtop}, we therefore add a new vertex $r$ and the arcs $(r,e)$ and $(r,v)$, and we put $\omega((r,v))=\frac{1}{2}(\lambda_5(\{v_5,e\})-l_{(N,\omega)}(v,c))=2$ and $\omega((r,e))=\frac{1}{2}(\lambda_5(\{v_5,e\})+l_{(N,\omega)}(v,c))=3$. This gives rise to the network depicted in Figure~\ref{fig:ex1}(ii), right, which is an ultrametric uprooting of the weighted phylogenetic tree depicted in Figure~\ref{fig:ex1}(i).

We now state the aforementioned result that any weighted phylogenetic
tree admits an ultrametric uprooting. As the proof is 
quite long and technical, we present it in the Appendix

\begin{proposition}\label{pr:uur}
Given a weighted phylogenetic tree $(T,\lambda)$ on $X$, the output $(N,\omega)$ of Algorithm \textsc{Ultrametric Uprooting} is an ultrametric uprooting $(N,\omega)$ of $(T,\lambda)$. In particular, all weighted phylogenetic trees admit an ultrametric uprooting.
\end{proposition}

We now turn our attention to the problem of proving that the ultrametric uprooting for a weighted phylogenetic tree $(T,\lambda)$ constructed by Algorithm~\ref{alg:ulup} is in fact the only possibly uprooting of $(T,\lambda)$ up to isomorphism. To do this, we first introduce some additional terminology.

Let $T$ be a phylogenetic tree with leaf set $X$, and let $N$ be an uprooting of $T$. As mentioned in the previous section, for all vertices $u,v$ of $T$, there is a unique path $P_T(u,v)$ in $T$ between $u$ and $v$. We say that a vertex $w$ of $P_T(u,v)$ distinct from $u$ and $v$ is a \emph{low point} of $P_T(u,v)$ (in $N$) if no proper descendant of $w$ in $N$ is a vertex of $P_T(u,v)$. In particular, a low point of $P_T(u,v)$ must have indegree $2$ or more in $N$. Indeed, since $w$ is distinct from $u$ and $v$, there exists two distinct vertices $w_1,w_2$ of $T$ that are adjacent to $w$ in $P_T(u,v)$. Since $w$ is a low point, neither $w_1$ nor $w_2$ is a child of $w$ in $N$. It follows that $w$ has two distinct parents $p_1,p_2$ in $N$, where for $i \in \{1,2\}$, $p_i$ is either $w_i$, or a root of $R_2(N)$ whose children are $w$ and $w_i$. Note however that not all vertices of $P_T(u,v)$ that have indegree $2$ or more in $N$ are low points. 
Moreover, $P_T(u,v)$ does not have any low point if and only if $x$ and $y$ share an ancestor in $N$.

Before proving our aforementioned uniqueness result, we present a useful lemma.
For $v$ a vertex of $T$, we define $cl_T(v)$ as the set of leaves of $T$ that are closest
to $v$ in $(T,\lambda)$ with regards to $\lambda$, that is, $cl_T(v) = cl_{(T,\lambda)}(v) =\mathrm{argmin}_{x \in X}\{l_{(T,\lambda)}(v,x)\}$.

\begin{lemma}\label{lm:clc}
Let $(T,\lambda)$ be a weighted phylogenetic tree on $X$, and let $(N,\omega)$ be an ultrametric uprooting of $(T,\lambda)$. Then for all vertices $v$ of $T$, we have $C_N(v)=cl_T(v)$.
\end{lemma}

\begin{proof}
To ease notation, we put $l=l_{(T,\lambda)}$ and $D=D_{(T,\lambda)}$. Recall that if $u,v$ are two vertices of $T$ such that $u$ is an ancestor of $v$ in $N$, then $l(u,v)=l_{(N,\omega)}(u,v)$.

We first remark that, since $(N,\omega)$ is ultrametric, $l(v,x)=l(v,x')$ holds for all $x,x' \in C_N(v)$. Hence, to show that $C_N(v)=\mathrm{argmin}_{x \in X}\{l_{(T,\lambda)}(v,x)\}$ holds, it suffices to show that for all pairs $x,y \in X$ such that $x \in C_N(v)$, $y \notin C_N(v)$, we have that $l(v,x)<l(v,y)$. So, let $x,y$ be such a pair. We show that $l(v,x)<l(v,y)$ holds by induction on the number of low points of $P_T(x,y)$.

To see the base case, suppose that $P_T(x,y)$ does not have any low point. Then, $x$ and $y$ share and ancestor in $N$. Let $w=\lca_N(x,y)$. Since $y \notin C_N(v)$, $w$ is not a descendant of $v$ in $N$.

If $w$ is an ancestor of $v$, we have $l(v,y)=l_{(N,\omega)}(w,y)+l_{(N,\omega)}(w,v)$. Moreover, $(N,\omega)$ is ultrametric, and so $l_{(N,\omega)}(w,y)=l_{(N,\omega)}(w,x)=l_{(N,\omega)}(w,v)+l_{(N,\omega)}(v,x)$. Combining these two equalities together, we obtain $l(v,y)=2l_{(N,\omega)}(w,v)+l_{(N,\omega)}(v,x)=2l(w,v)+l(v,x)$. Hence, $l(v,y)>l(v,x)$ holds in this case.

If otherwise, $w$ is not an ancestor of $v$, we denote by $h$ the first vertex of $N$ that is common to the paths from $v$ to $x$ and from $w$ to $x$ in $N$. Then we have $l(v,y)=l_{(N,\omega)}(w,y)+l_{(N,\omega)}(w,h)+l_{(N,\omega)}(v,h)$. Moreover, $(N,\omega)$ is ultrametric, so $l_{(N,\omega)}(w,y)=l_{(N,\omega)}(w,x)=l_{(N,\omega)}(w,h)+l_{(N,\omega)}(h,x)$. Putting these two equalities together, we obtain $l(v,y)=2l_{(N,\omega)}(w,h)+l_{(N,\omega)}(v,h)+l_{(N,\omega)}(h,x)=2l(w,h)+l(v,x)$. Hence, $l(v,y)>l(v,x)$ also holds in this case.

Now, suppose that $P_T(x,y)$ contains $h \geq 1$ low points, and that 
for all $y' \in X$ such that $y' \notin C_N(v)$ and $P_T(x,y')$ has $h'<h$ low points, we have that $l(v,x)<l(v,y')$.
In addition, let $u$ be a low point of $P_T(x,y)$ such that $u$ and $y$ 
share an ancestor in $N$, and let $z$ be 
a descendant of $u$ in $N$. 

Clearly, $u$ is a vertex of both $P_T(z,v)$ and $P_T(y,v)$. In particular, we have $l(v,y)=l(y,u)+l(u,v)=l(y,z)+l(z,v)-2l(u,z)$. Moreover, $l(y,z)=\D(y,z)=2l_{(N,\omega)}(w,z)$, where $w=\lca(y,z)$, so the previous equality can be written $l(v,y)=2l_{(N,\omega)}(w,z)+l(z,v)-2l_{(N,\omega)}(u,z)$. Since $u$ is not an ancestor of $y$, but shares an ancestor with $y$, it follows that $w$ is a proper ancestor of $u$ in $N$. Hence, $l_{(N,\omega)}(w,z)>l_{(N,\omega)}(u,z)$. It follows that $l(v,y)>l(v,z)$. By choice of $z$, the path $P_T(x,z)$ has $h-1$ low points. By our induction hypothesis, $l(v,z)>l(v,x)$ follows, so $l(v,y)>l(v,x)$ holds as desired.
\end{proof}

We are now ready to prove our uniqueness result.

\begin{proposition}\label{pr:uniqueness}
Let $(T,\lambda)$ be a weighted phylogenetic tree on $X$. Then up to isomorphism, there exists a unique ultrametric uprooting of $(T,\lambda)$. 
\end{proposition}

\begin{proof}
Let $(N,\omega)$ be an ultrametric uprooting of $(T,\lambda)$. For $\{u,v\}$ an edge of $T$, exactly one of the following must hold: $(u,v)$ is an arc of $N$, $(v,u)$ is an arc of $N$, or there exists a root $r \in R_2(N)$ such that $u$ and $v$ are the children of $r$ in $N$. To show that $N$ is uniquely determined by $(T,\lambda)$, we show that for all edges of $N$, the choice between the aforementioned three possibilities is uniquely determined by $(T,\lambda)$.

So, let $\{u,v\}$ be an edge of $T$. By Lemma~\ref{lm:clc}, all leaves $x \in cl_T(v)$ are descendants of $v$ in $N$. In particular, if there exists some $x \in cl_T(v)$ such that $u$ is a vertex of $P_T(v,x)$, then $(v,u)$ must be an arc in $N$. By symmetry, if there exists $y \in cl_T(u)$ such that $v$ is a vertex of $P_T(u,y)$, then $(u,v)$ must be an arc in $N$. Note that these two situations cannot happen simultaneously. Indeed, suppose for contradiction that there exist $x$ and $y$ as specified. Then, we have $l_{(T,\lambda)}(v,x)=l_{(T,\lambda)}(u,x)+\lambda(\{u,v\})>l_{(T,\lambda)}(u,x)$ and $l_{(T,\lambda)}(u,y)=l_{(T,\lambda)}(v,y)+\lambda(\{u,v\})>l_{(T,\lambda)}(v,y)$. Since $x \in cl_T(v)$ and $y \in cl_T(u)$, we also have $l_{(T,\lambda)}(v,x) \leq l_{(T,\lambda)}(v,y)$ and $l_{(T,\lambda)}(u,y) \leq l_{(T,\lambda)}(u,x)$. Taken together, these four inequalities yield a contradiction.

We next show that if there is no $x \in cl_T(v)$ such that $u$ is a vertex of $P_T(v,x)$ and no $y \in cl_T(u)$ such that $v$ is a vertex of $P_T(u,y)$, then there exists $r \in R_2(N)$ such that $u$ and $v$ are the children of $r$ in $N$. To see this, suppose for contradiction that one of $(u,v)$ or $(v,u)$, say $(v,u)$, is an arc of $N$. Then, we have $C_N(u) \subsetneq C_N(v)$. So, let $x \in C_N(u)$. By Lemma~\ref{lm:clc}, we have $x \in cl_T(v)$. Moreover, $x \in C_N(u)$ and $u$ is a child of $v$, so the path from $v$ to $x$ in $N$ contains $u$. Hence, there is a path in $T$ between $v$ and $x$ that contains $u$. Such a path being unique, it follows that $u$ is a vertex of $P_T(v,x)$, a contradiction to our assumption that there is no leaf $x \in cl_T(v)$ such that $u$ is a vertex of $P_T(v,x)$.

We have shown that up to isomorphism, $N$ is uniquely determined by $(T,\lambda)$. It remains to show that the weighting $\omega$ is also uniquely determined by $(T,\lambda)$. Since $(N,\omega)$ is a weight-preserving uprooting of $(T,\lambda)$, $\omega((u,v))=\lambda(\{u,v\})$ holds by definition for all arcs $(u,v)$ of $N$ such that $u \notin R_2(N)$. Consider a vertex $r \in R_2(N)$ with children $u$ and $v$. By definition, $\{u,v\}$ is an edge of $T$, and we have $\omega((r,u))+\omega((r,v))=\lambda(\{u,v\})$. Moreover, $(N,\omega)$ is an ultrametric arboreal network, so for all $x \in C_N(u)$, $y \in C_N(v)$, $l_{(N,\omega)}(r,x)=l_{(N,\omega)}(r,y)$ holds. Since $l_{(N,\omega)}(r,x)=\omega((r,u))+l_{(N,\omega)}(u,x)$ and $l_{(N,\omega)}(r,y)=\omega((r,v))+l_{(N,\omega)}(v,y)$, we have $\omega((r,u))+l_{(N,\omega)}(u,x)=\omega((r,v))+l_{(N,\omega)}(v,y)$. Together with the previous equality, it follows that $\omega((r,u))=\frac{1}{2}(\lambda(\{u,v\})+l_{(N,\omega)}(v,y)-l_{(N,\omega)}(u,x))$ and $\omega((r,v))=\frac{1}{2}(\lambda(\{u,v\})+l_{(N,\omega)}(u,x)-l_{(N,\omega)}(v,y))$. Since the length of the paths from $u$ to $x$ and from $v$ to $y$ are uniquely determined by $N$ and $(T,\lambda)$, this is also the case of $l_{(N,\omega)}(u,x)$ and $l_{(N,\omega)}(v,y)$. Hence, $\omega((r,u))$ and $\omega((r,v))$ are uniquely determined by $(T,\lambda)$. This concludes the proof that $\omega(e)$ is uniquely determined by $(T,\lambda)$ for all arcs $e$ of $N$.
\end{proof}

An important consequence of Proposition~\ref{pr:uniqueness} is that 
	choices
	made by algorithm \textsc{Ultrametric Uprooting}
	do not affect its output.
	More specifically, we have:

\begin{observation}
The output of Algorithm~\ref{alg:ulup} is independent of the choice of the wcps at Line~\ref{l:wcps}, and of the successive choices of the element $x$ at Line~\ref{l:defvx}.
\end{observation}

Taking Propositions~\ref{pr:uur} and \ref{pr:uniqueness} together with Theorem~\ref{th:treelike}, we immediately 
obtain the main result of this section:

\begin{theorem}\label{th:rooting}
Let $D$ be a distance on $X$. The following are equivalent:
\begin{itemize}
\item[(i)] There exists an ultrametric arboreal network $(N,\omega)$ on $X$ such that $D=D_{(N,\omega)}$.
\item[(ii)] There exists a weighted phylogenetic tree $(T,\lambda)$ on $X$ such that $D=D_{(T,\lambda)}$.
\item[(iii)] For all $x,y,z,u \in X$, 
$$D(x,y)+D(z,u) \leq \mathrm{max}\{D(x,z)+D(y,u),D(x,u)+D(y,z)\}.$$
\end{itemize}
Moreover, if this holds, then both $(T,\lambda)$ and $(N,\omega)$ are unique up to isomorphism, and $(N,\omega)$ is an ultrametric uprooting of $(T,\lambda)$.
\end{theorem}

\section{Characterizing arboreal ultrametrics}\label{sec:pd}

In this section we turn our attention to giving a 
characterization for when a partial distance is an arboreal ultrametric. We begin with some definitions.

A map $\D: X \times X \to \mathbb R_{\geq 0} \cup \{\infty\}$ is called
a {\em partial distance (on $X$)}  if it is symmetric and vanishes precisely on the diagonal.
To any weighted arboreal network $(N,\omega)$ on $X$, we can associate a partial distance $\D_{(N,\omega)}$ on $X$ as follows: if $x,y \in X$ share an ancestor in $N$, we put $\D_{(N,\omega)}(x,y)=l_{(N,\omega)}(v,x)+l_{(N,\omega)}(v,y)$, where $v=\lca_N(x,y)$. Note that in this case $\D_{(N,\omega)}(x,y)=D_{(N,\omega)}(x,y)$ and, if $(N,\omega)$ is ultrametric, then 
$l_{(N,\omega)}(v,x)=l_{(N,\omega)}(v,y)$, and so $\D_{(N,\omega)}(x,y)=2l_{(N,\omega)}(v,x)=2l_{(N,\omega)}(v,y)$ also holds. If otherwise, $x$ and $y$ do not share an ancestor in $N$, we put $\D_{(N,\omega)}(x,y)=\infty$. Note that $\D_{(N,\omega)}$ is a distance on $X$ if and only if $N$ is a rooted phylogenetic tree, in which case $\D_{(N,\omega)}$ is
equal to $D_{(N,\omega)}$. We say that a partial distance $\D$ on $X$ is \emph{arboreal representable} if there exists an ultrametric arboreal network $(N,\omega)$ on $X$ satisfying $\D_{(N,\omega)}=\D$. In this case, we say that $(N,\omega)$ \emph{represents} $\D$.

Partial distances induced by weighted arboreal networks are closely related to the symbolic arboreal maps introduced in \cite{huber2024shared}. For $M$ be a non-empty set of {\em symbols} and $\odot$ an element that is not in $M$,
a \emph{symbolic map} on $X$ is a map
$d: {X \choose 2} \to M \cup \{\odot\}$. In particular, we shall view a
partial distance $\D$ on $X$ as a symbolic ultrametric $d$ where $M=\mathbb R_{>0}$ and $\odot = \infty$,
where the fact that $\D$ is symmetric ensures that $d$ is well defined.


In what follows, we will make use of a characterization of 
symbolic maps that arise from arboreal networks given in \cite{huber2024shared}, which we now recall for the 
convenience of the reader.
Given an arboreal network $N$ on $X$, we denote by $V(N)^-$ the set of all vertices of $N$ of outdegree $2$ or more. A \emph{labelled arboreal network} is a pair $(N,t)$ where $N$ is an arboreal network and $t: V(N)^- \to M$ is a map assigning an element of $M$ to each vertex of $V(N)^-$. For a symbolic map $d: {X \choose 2} \to M \cup \{\odot\}$,
define the undirected graph $G_d$ to be the graph with vertex set $X$ and edges precisely those 
$\{x,y\} \in {X \choose 2}$ such that $d(x,y) \neq \odot$.
Then, we say that $(N,t)$ \emph{explains} $d$ if, for all $x,y \in X$ distinct, $d(x,y)=t(\lca_N(x,y))$ if $x$ and $y$ share an ancestor in $N$, and $d(x,y)=\odot$ otherwise. Symbolic maps that can be explained by labelled arboreal networks, also called as \emph{symbolic arboreal maps}, were characterised in \cite{huber2024shared} as follows:

\begin{theorem}[\cite{huber2024shared}, Theorem 7.5]\label{thm:diss}
Suppose that $X$ is a set with $|X|\geq 2$ and that $d: {X \choose 2} \to M \cup \{\odot\}$ is a symbolic map. Then, 
$d$ is a symbolic arboreal map if and only if the following four properties all hold:
\begin{itemize}
\item[(A1)] $G_d$ is connected and Ptolemaic.
\item[(A2)] No three elements $x,y,z \in X$ satisfy $|\{d(x,y),d(x,z),d(y,z)\}|=3$ and\\ $\odot \notin \{d(x,y),d(x,z),d(y,z)\}$.
\item[(A3)] No four elements $x,y,z,u \in X$ satisfy $d(x,y)=d(y,z)=d(z,y) \neq d(y,u)=d(u,x)=d(x,z)$ and $\odot \notin \{d(x,y),d(x,z)\}$.
\item[(A4)] For all $x,y,z,u \in X$ distinct such that $d(z,u)= \odot$ and $d$ maps all other elements of ${\{x,y,z,u\}\choose 2}$ to an element of $M$, both $d(x,z)=d(y,z)$ and $d(x,u)=d(y,u)$ hold.
\end{itemize}
\end{theorem}

We now use this result to characterize arboreal-representable partial distances.  First, we show that if $(N,\omega)$ is an ultrametric arboreal network, then 
$\D_{(N,\omega)}$ 
is a symbolic arboreal map.

\begin{lemma}~\label{le:ds}
Let $(N,\omega)$ be an ultrametric arboreal network on $X$. Then, there exists a labelling map $t:V(N)^-\to \mathbb R_{>0}$ such that 
	the labelled network $(N,t)$ explains $\D_{(N,\omega)}$.  In particular, $\D_{(N,\omega)}$ is a symbolic arboreal map on $X$.
\end{lemma}

\begin{proof}
To ease notation, we put $\D=\D_{(N,\omega)}$. As already mentioned, $\D$ is a symbolic map on $X$. To obtain the labelling map $t$, 
	we put, for all $v\in V(N)^-$,
 $t(v)=2l_{(N,\omega)}(v,x)$, where $x \in C_N(v)$. Note that, 
since $(N,\omega)$ is ultrametric, $t(v)$ does not depend on the choice of $x$ in $C_N(v)$. 

We now show that $(N,t)$ explains $\D$.
Let $x,y \in X$ distinct. Since $(N,\omega)$ represents $\D$, we 
have $\D(x,y)= \infty$ if and only if $x$ and $y$ do not share an ancestor. Suppose now that $\D(x,y) \not= \infty$. Then $x$ and $y$ have a common ancestor in $N$. Let $v=\lca_N(x,y)$. By definition, $\D(x,y)=l_{(N,\omega)}(v,x)+l_{(N,\omega)}(v,y)$. Since $(N,\omega)$ is ultrametric, we have $l_{(N,\omega)}(v,x)=l_{(N,\omega)}(v,y)$, so $D_{(N,\omega)}(x,y)=2l_{(N,\omega)}(v,x)=t(v)$. Hence, $(N,t)$ explains $\D$.
\end{proof}

We next present a characterization for arboreal representable partial distances.

\begin{theorem}\label{th:ultra}
Let $\D: X \times X \to \mathbb R_{>0} \cup \{\infty\}$ be a partial distance on $X$. Then there exists an ultrametric 
arboreal network on $X$ representing $\D$ if and only if the following 
three properties hold:
\begin{itemize}
\item[(U1)] $G_{\D}$ is connected and chordal.
\item[(U2)] If $x,y,z\in X$ are pairwise distinct such that $\infty \notin \{\D(x,y),\D(x,z),\D(y,z)\}$, then $\D(x,y) \leq \mathrm{max}\{\D(x,z),\D(y,z)\}$.
\item[(U3)] If $x,y,z,u\in X$ are pairwise distinct such that $\D(z,u)= \infty$ 
and $\D$ maps all 
other elements of ${\{x,y,z,u\}\choose 2}$ to an element in $\mathbb R_{>0}$, then 
$$\D(x,y)< \mathrm{min}\{\D(x,z), \D(x,u), \D(y,z), \D(y,u)\}.$$
\end{itemize}
\end{theorem}

The intuition behind Property~(U3) is as follows. If $N$ is an ultrametric arboreal network on $X$ and $x$, $y$, $z$ and $u$ are as in Property~(U3) then if $x$ and $y$ share an ancestor with $z$ and an ancestor with $u$, but $z$ and $u$ do not share an ancestor, then there must exist a vertex of indegree $2$ or more in $N$ that is an ancestor of $x$ and $y$ but not of $z$ and $u$. As a consequence, $x$ and $y$ are "closer`` 
	to each other than they are to $z$ and to $u$, respectively.

To further illustrate Property~(U3), consider the distance $\D$ on $X=\{x,y,z,u\}$ defined by $\D(x,z)=1$, $\D(y,z)=\D(x,y)=2$, $\D(x,u)=\D(y,u)=3$, and $\D(z,u)=\infty$. Then it is straight-forward to verify that $\D$ satisfies (U1) and (U2), but not (U3). If there existed an ultrametric arboreal network $N$ on $X$ representing $\D$, then the restriction of $\D$ to $\{x,y,z\}$ implies that the last common ancestor of $x$ and $y$ is an ancestor of the last common ancestor of $x$ and $z$, and the restriction of $\D$ to $\{x,y,u\}$ implies that the last common ancestor of $x$ and $u$ is an ancestor of the last common ancestor of $x$ and $y$. Taken together, these two observations imply that $z$ and $u$ share an ancestor in $N$ contradicting $\D(z,u)=\infty$.

Before proving Theorem~\ref{th:ultra}, we make an observation about Property~(U3) that will be useful in the proof and in subsequent results.

\begin{observation}\label{obs:gem}
If $\D$ satisfies (U3), then $G_{\D}$ is gem-free and $W_5$-free. In particular, if $\D$ satisfies (U1) and (U3), then $G_{\D}$ is Ptolemaic. \end{observation}

\begin{proof}
Suppose that $\D$ satisfies (U3), and assume for contradiction that there exists pairwise distinct vertices $x,y,z,u,v$ of $G_{\D}$ such that the subgraph of $G_{\D}$ induced by these five vertices is a gem or a $W_5$. Then up to a permutation, $G_{\D}$ contains the edges $\{x,y\},\{x,v\},\{y,z\},\{y,v\};\{z,u\},\{z,v\}$ and $\{u,v\}$, and does not contain the edges $\{x,z\}$ and $\{y,u\}$. Applying (U3) on $x,y,z,u$, we obtain $\D(y,v)<\D(z,v)$, and applying (U3) on $y,z,u,v$, we get $\D(z,v)<\D(y,v)$. This is impossible, so $G_{\D}$ is gem-free and $W_5$-free. If in addition, $\D$ satisfies (U1), then $G_{\D}$ is chordal, so $G_{\D}$ is Ptolemaic.
\end{proof}

We now proceed with the proof of Theorem~\ref{th:ultra}.

\begin{proof}
Suppose first that there exists an ultrametric arboreal network $(N,\omega)$ on $X$ representing $\D$. By Lemma~\ref{le:ds}, $\D$ is a symbolic arboreal map, and so it follows by Theorem~\ref{thm:diss} that $\D$ satisfies Properties~(A1)-(A4). Since a Ptolemaic graph is chordal, (U1) is a weaker version of (A1), so Property~(U1) holds.

To see that Property~(U2) holds, let $x,y,z \in X$ be three pairwise distinct elements such that $\infty \notin \{\D(x,y),\D(x,z),\D(y,z)\}$. Since $(N,\omega)$ represents $\D$ and is arboreal, it follows that there exists a root $r$ of $N$ such that $x,y,z \in C_N(r)$. In particular, the restriction of $\D$ to $C_N(r)$ is an ultrametric on $C_N(r)$. By Theorem~\ref{th:ultree}, $\D(x,y) \leq \mathrm{max}\{\D(x,z),\D(y,z)\}$ follows.

To see that Property~(U3) holds, let $x,y,z,u \in X$ be four pairwise distinct elements such that $\D(z,u)= \infty$ and $\D$ 
maps all other elements in ${\{x,y,z,u\}\choose 2}$ to an element in $\mathbb R_{>0}$. Since $\D$ satisfies Property~(A4), we have $\D(x,z)=\D(y,z)$ and $\D(x,u)=\D(y,u)$. Moreover, since $\D$ satisfies Property~(U2) by the previous paragraph, we have $\D(x,z)=\D(y,z) \geq \D(x,y)$ and $\D(x,u)=\D(y,u) \geq \D(x,y)$.
It remains to show that these two inequalities are strict. 

To see that this is the case, suppose for contradiction that $\D(x,z)=\D(x,y)=\D(y,z)$. Since 
$N$ is arboreal, and $\omega(a)>0$ for all arcs $a$ of $N$, $\lca_N(x,y)=\lca_N(x,z)=\lca_N(y,z)$ must hold. Denote this least common ancestor by $v$. 
Since $\D(x,u), \D(y,u) \in \mathbb R_{>0}$, it follows 
that $x, y$ and $u$ share an ancestor in $N$. Then, either $u$ is a descendant of $v$ in $N$, or $v$ is a descendant of $\lca_N(x,u)$. In the first case, $v$ is a common ancestor of $z$ and $u$, and in the second case, $\lca_N(x,u)$ is a common ancestor of $z$ and $u$. Both are impossible, since $\D(z,u)=\infty$ implies that $z$ and $u$ do not share an ancestor in $N$. Hence, $\D(x,z)=\D(y,z) > \D(x,y)$ must hold. By symmetry, we also have that $\D(x,u)=\D(y,u) > \D(x,y)$ holds.

Conversely, suppose 
that $\D$ satisfies Properties~(U1), (U2) and (U3). We begin by showing 
that $\D$ satisfies Properties~(A1) to (A4).
Since $\D$ satisfies Property~(U1), $G_{\D}$ is connected and chordal. Moreover, by Observation~\ref{obs:gem}, $G_{\D}$ is gem-free. Hence, $G_{\D}$ is Ptolemaic, so Property~(A1) holds.

We now show that Property~(A2) is a consequence of Property~(U2). Let $x,y,z \in X$ be such that $\odot=\infty \notin \{\D(x,y),\D(x,z),\D(y,z)\}$. Up to permutation of the elements in $\{x,y,z\}$, we may assume that $\D(x,y)=\mathrm{max}\{\D(x,y),\D(x,z),\D(y,z)\}$. By (U2), $\D(x,y) \leq \mathrm{max}\{\D(x,z),\D(y,z)\}$, so by choice of the pair $x,y$, it follows that $\D(x,y) = \mathrm{max}\{\D(x,z),\D(y,z)\}$ follows. Hence, we have $|\{\D(x,y),\D(x,z),\D(y,z)\}| < 3$, so (A2) holds.

Next, we show that Property~(A3) holds. Let $x,y,z,u \in X$ be such that $\D(x,y)=\D(y,z)=\D(z,u) \neq \D(z,x)=\D(x,u)=\D(u,y)$ and $\infty \notin \{\D(x,y),\D(x,z)\}$. 
Using Property~(U2) on the sets $\{x,y,z\}$ and $\{x,y,u\}$, we must have $\D(x,y)=\D(y,z)>D(x,z)$ and $\D(x,u)=\D(u,y)>\D(x,y)$, which is impossible 
since $\D(x,z)=\D(x,u)$. Hence, Property~(A3) holds.

Finally, we show that $\D$ satisfies Property~(A4). Let $x,y,z,u \in X$ be pairwise distinct such that $d(z,u)=\odot$ while $d$ maps all other pairs of elements of ${\{x,y,z,u\}}$ to an element of $M=\mathbb R_{>0}$. By (U3), we have $\D(x,y)< \mathrm{min}\{\D(x,z), \D(x,u), \D(y,z), \D(y,u)\}$. Moreover, (U2) applied to the set $\{x,y,z\}$ gives $\D(y,z) \leq \mathrm{max}\{\D(x,y),\D(x,z)\}=\D(x,z)$ and $\D(x,z)  \leq \mathrm{max}\{\D(x,y),\D(y,z)\}=\D(y,z)$, so $\D(x,y)<\D(x,z)=\D(y,z)$ follows. Applying (U2) to the set $\{x,y,u\}$ implies $\D(x,y)<\D(x,u)=\D(y,u)$ in a similar fashion.

Since $\D$ satisfies Properties~(A1)-(A4), Theorem~\ref{thm:diss} implies 
that there exists a labelled arboreal network $(N,t)$ that explains $\D$. 
Without loss of generality, we may assume that $(N,t)$ is such that no arc $(u,v)$ of $N$ 
satisfies $t(u)=t(v)$ as otherwise we could contract that arc and the resulting arboreal network would also explain $\D$. We next construct a weighting $\omega$ of $N$ such that $(N,\omega)$ represents $\D$.
To this end, we first define a map $\delta: V(N) \to \mathbb R_{\geq 0}$ as follows. If $v$ is a leaf of $N$, 
we put $\delta(v)=0$. If $v$ has outdegree at least $2$, we put $\delta(v)=\frac{1}{2}t(v)$. 
To be able to extend the definition of $\delta$ to $V(N)$, we claim
that for all vertices $u$ and $v$ of $N$ distinct such that both $u$ and $v$ have outdegree at least 2 
and $u$ is an ancestor of $v$ in $N$, we have $\delta(u)>\delta(v)$. 

Let $u$ and $v$ be two such vertices of $N$. First, note that if $u=v_1, \ldots, v_k=v$, $k \geq 2$, is the subsequence of vertices of outdegree at least $2$ 
on the path $P_N(u,v)$ from $u$ to $v$, then it suffices to show that $\delta(v_i) > \delta(v_{i+1})$ holds for all $1 \le i \le k-1$.
By definition of the vertices $v_1, \ldots, v_k$, there is no vertex of outdegree $2$ or more on the path from $v_i$ to $v_{i+1}$. 
Hence to prove the claim, we may assume without loss of generality that $u$ and $v$ 
are such that other than possibly $u$ or $v$, no other vertex on $P_N(u,v)$ 
has outdegree 2 or more in $N$.

Now, suppose $x,y  \in C_N(v)$  are such that $v=\lca_N(x,y)$, and  let $z \in C_N(u)$ be such that $z \notin C_N(v)$. 
Note that since both $u$ and $v$ have outdegree at least $2$ in $N$, and $N$ is arboreal, the leaves $x,y$ and $z$ always exist. 
In particular, $u=\lca_N(x,z)=\lca_N(y,z)$. 
Since $(N,t)$ explains $\D$, it follows that $\D(x,y)=t(v)$ and $\D(x,z)=\D(y,z)=t(u)$. By 
Property~(U2), $t(u) \geq t(v)$ holds. We now show that this inequality is strict
which immediately implies our claim by the definition of the map $\delta$. To do this we distinguish between two cases: 
(a) $(u,v)$ is an arc of $N$, and (b) there exists a vertex $h$ of indegree $2$ or more on the path from $u$ to $v$.

If case (a) holds, we have $t(u) \neq t(v)$ by assumption on $(N,t)$, so $t(u) > t(v)$. 
Hence, $\delta(u)>\delta(v)$ holds as claimed. If case (b) holds, let $r$ be a root of $N$ 
such that $h$ is a descendant of $r$ and $u$ is not, and let $z' \in X$ be a descendant of $r$ 
that is not also a descendant of $h$. Note that $x,y,z,z'$ are pairwise distinct. Since $(N,t)$ explains $\D$, we have $\D(z,z')=\infty$, while 
$\D$ maps all other elements in ${\{x,y,z,z'\}\choose 2}$ to an element in $\mathbb R_{>0}$. By Property~(U3), it follows 
that $t(u)=\D(x,z)=\D(y,z)>\D(x,y)=t(v)$, so $\delta(u)>\delta(v)$ also holds in this case. This completes the proof of the claim.

The claim being true, it is therefore always possible to extend the definition of $\delta$ to $V(N)$  in such a way that $\delta(u)>\delta(v)$ holds
for all arcs $(u,v)$ of $N$. Putting $\omega(a)=\delta(u)-\delta(v)$ for all arcs $a=(u,v)$ of $N$, it follows that $\omega(a)>0$. Thus, $(N,\omega)$ is a weighted arboreal network. 

We now show that $(N,\omega)$ is an ultrametric arboreal network. Let 
$v$ be a non-leaf vertex of $N$ and let $x \in X$ such that $x$ is a descendant of $v$ in $N$. 
Let $v_1=v,v_2, \ldots, v_k=x$, $k \geq 2$, be a path from $v$ to $x$ in $N$. By definition, we have $l_{(N,\omega)}(v,x)=\sum_{i=1}^{k-1} \omega((v_i,v_{i+1}))$. 
Since $\omega((v_i,v_{i+1}))=\delta(v_i)-\delta(v_{i+1})$, and $\delta(x)=0$, 
it follows that $l_{(N,\omega)}(v,x)=\delta(v)$. In particular, 
$l_{(N,\omega)}(v,x)$ does not depend on the choice of $x$ in $C_N(v)$, so $(N,\omega)$ is an ultrametric arboreal network.

To complete the proof, it remains to show that $(N,\omega)$ represents $\D$. 
Since $\D$ is a partial distance  and $(N,t)$ explains $\D$, we have $\D(x,y)=\infty$ if and only if $x$ and $y$ do not 
share an ancestor in $N$. Now, let $x,y \in X$ distinct be such that $\D(x,y) \neq \infty$, and let $v=\lca_N(x,y)$. By definition, and since $(N,\omega)$ is ultrametric, we have $D_{(N,\omega)}(x,y)=l_{(N,\omega)}(v,x)+l_{(N,\omega)}(v,y)=2l_{(N,\omega)}(v,x)$. As 
observed in the previous paragraph, $l_{(N,\omega)}(v,x)=\delta(v)$, and $\delta(v)=\frac{1}{2}t(v)$ by definition 
of $\delta$. Since $t(v)=\D(x,y)$ as $(N,t)$ explains $\D$ and $v=\lca_N(x,y)$, it follows that $D_{(N,\omega)}(x,y)=\D(x,y)$. This concludes the proof of the theorem.
\end{proof}

A first immediate consequence of Theorem~\ref{th:ultra} is the following, computational complexity result:

\begin{corollary}\label{cor:on4}
For $\D$ a partial distance on $X$, we can decide in time $O(|X|^4)$ whether or not $\D$ is an arboreal ultrametric.
\end{corollary}

The following corollary characterizes when 
the restriction of an arboreal ultrametric on $X$ to a subset of $X$ is also an arboreal ultrametric.

\begin{corollary}\label{cor:hered}
Let $\D$ be an arboreal ultrametric on $X$, and let $Y$ be a nonempty subset of $X$. The restriction $\D_Y$ of $\D$ to $Y \times Y$ is an arboreal ultrametric if and only if $G_{\D_Y}$ is connected.
\end{corollary}

\begin{proof}
Since $\D$ is an arboreal ultrametric, $\D$ satisfies Properties (U1), (U2) and (U3) of Theorem~\ref{th:ultra}.

Since $\D$ satisfies (U1), $G_{\D}$ is connected and chordal. Moreover, $G_{\D_Y}$ is precisely the subgraph of $G_{\D}$ induced by the elements of $Y$. In particular, $G_{\D_Y}$ is chordal. Moreover, it is straightforward to verify that if $\D$ satisfies (U2) (\emph{resp.} (U3)), then $\D_Y$ also satisfies (U2) (\emph{resp.} (U3)). In summary, $\D_Y$ satisfies (U2) and (U3), and $G_{\D_Y}$ is chordal. By Theorem~\ref{th:ultra}, $\D$ is an arboreal ultrametric if and only if $G_{\D_Y}$ is connected.
\end{proof}

We conclude this section with a uniqueness result. For the purpose of this result, we will relax the definition of a network, by allowing leaves to have degree more than one. Note that if a network $N$ contains such leaves, it can be transformed into a network in the 
sense given above by applying the following three steps for all leaves $x$ of indegree $2$ or more. First, introduce a new vertex $v$. Then, replace all arcs $(u,x)$ with the arc $(u,v)$. Finally, add the arc $(v,x)$. Note that the network $N^+$ obtained this way is unique. However, if $(N,\omega)$ is edge weighted, then there are infinitely many ways to assign weights $\omega'$ to the arcs of $N^+$ in such a way that $(N^+,\omega')$ is a weighted network satisfying $\D_{(N^+,\omega')}=\D_{(N,\omega)}$.

To state our uniqueness result, we shall use the following fact which follows by \cite[Theorem 7.6]{huber2024shared}.
Suppose that $d$ is a symbolic map, then there exists a unique (up to isomorphism) labelled arboreal network $(N,t)$ such that $(N,t)$ explains $d$, 
$N$ does not contain vertices of outdegree $1$, and 
$t(u) \neq t(v)$ for all arcs $(u,v)$ of $N$ with $u,v \in V(N)^-$.

\begin{theorem}\label{th:uniq}
Let $\D: X \times X \to \mathbb R_{>0} \cup \{\infty\}$ be a partial distance on $X$. If $\D$ satisfies properties (U1) to (U3), then, up to isomorphism, there exists a unique ultrametric arboreal 
 network $(N, \omega)$ representing $\D$ such that $N$ does not contain vertices of outdegree $1$.
\end{theorem}

\begin{proof}
Since $\D$ satisfies properties (U1) to (U3), Theorem~\ref{th:ultra} implies that there exists an ultrametric arboreal network representing $\D$. Let $(N,\omega)$ be such a network, and suppose that $N$ contains at least one vertex with outdegree $1$. Let $u$ be one such vertex, $v$ denote the unique child of $u$, and let $p_1, \ldots, p_k$, $k \ge 2$, be the parents of $u$ (note that $k \geq 2$ since $u$ has outdegree $1$ and $N$ is a network, and therefore $u$ is not a root and it does not have indegree $1$). Consider the weighted network $(N',\omega')$ where $N'$ is obtained by contracting the arc $(u,v)$, and where $\omega'(a)$ is defined, for all arcs $a$ of $N'$ as $\omega'(a)=\omega(a)+\omega((u,v))$ if $a=(p_i,u)$ for some $i \in \{1, \ldots, k\}$, and $\omega'(a)=\omega(a)$. It is straightforward to check that $(N',\omega')$ remains an ultrametric arboreal network explaining $\D$. By repeatedly applying this operation to vertices 
with outdegree $1$, we can therefore assume without loss of generality that  $(N,\omega)$ does not contain vertices of outdegree $1$.

Now, let $t_{\omega}: V(N)^- \to \mathbb R_{>0}$ be the map defined, for all vertices $v$ of $N$ of outdegree $2$ or more (that is, for all non-leaf vertices of $N$), by $t_{\omega}(v)=2l_{(N,\omega)}(v,x)$, where $x \in C_N(v)$. Recall that since $N$ is ultrametric, $l_{(N,\omega)}(v,x)$ does not depend on the choice of $x$ in $C_N(v)$, so $t_{\omega}(v)$ is uniquely determined.

Clearly, $(N,t_{\omega})$ explains the symbolic map $\D$. Indeed, since $(N,\omega)$ represents $\D$ as a partial distance, $\D(x,y)=\infty$ for some $x,y \in X$ if and only if $x$ and $y$ do not share an ancestor in $N$. Otherwise, if $x$ and $y$ do share an ancestor in $N$, then for $v=\lca_N(x,y)$, we have $\D(x,y)=l_{(N,\omega)}(v,x)+l_{(N,\omega)}(v,y)$, and since $(N,\omega)$ is ultrametric, $l_{(N,\omega)}(v,x)=l_{(N,\omega)}(v,y)$ and hence $\D(x,y)=2l_{(N,\omega)}(v,x)=t_{\omega}(v)$.

Now, by the fact mentioned before the statement of the 
theorem, there exists a unique (up to isomorphism) labelled arboreal network $(N_0,t_0)$ such that $(N_0,t_0)$ explains $\D$, $N_0$ 
does not contain vertices of outdegree $1$, and 
$t_0(u) \neq t_0(v)$ for all arcs $(u,v)$ of $N_0$ with $u,v \in V(N)^-$.
But as has already been established, $(N,t_{\omega})$ satisfies (i) and (ii). Moreover, $(N,t_{\omega})$ also satisfies (iii). Indeed, if $(u,v)$ is an arc of $N$ such that $u,v \in V(N)^-$, then $t_{\omega}(u)=t_{\omega}(v)+\omega((u,v))$ and, since $\omega((u,v))>0$, $t_{\omega}(u) \neq t_{\omega}(v)$ follows.
Hence $N$ and $N_0$ must be isomorphic.

In light of this last observation, it follows that if
$(N',\omega')$ is an ultrametric arboreal network representing $\D$ 
that does not contain any vertices of outdegree $1$,
then $N'$ is isomorphic to  $N$.
It remains to show that the weighting $\omega$ is uniquely determined by $N$ and $\D$. So, let $(u,v)$ be an arc of $N$. Suppose first that $v$ is a leaf of $N$. Since $u$ has outdegree $2$ or more, there exists a leaf $x \in C_N(u)$ distinct from $v$. Since $(N,\omega)$ represents $\D$, we have $\D(x,v)=l_{(N,\omega)}(u,x)+l_{(N,\omega)}(u,v)$. Moreover, $(N,\omega)$ is ultrametric, so we have $l_{(N,\omega)}(u,x)=l_{(N,\omega)}(u,v)$. Since $l_{(N,\omega)}(u,v)=\omega((u,v))$, it follows that $\omega((u,v))=\frac{1}{2}\D(x,v)$.

Suppose now that $v$ is not a leaf of $N$. Let $x,y \in C_N(v)$ distinct such that $v=\lca_N(x,y)$, and let $z \in C_N(u) \setminus C_N(v)$. Note that the existence of $x$, $y$ and $z$ is guaranteed by the fact that both $u$ and $v$ have outdegree $2$ or more in $N$. Since $(N,\omega)$ represents $\D$, we have $\D(x,y)=l_{(N,\omega)}(v,x)+l_{(N,\omega)}(v,y)$ and $\D(x,z)=l_{(N,\omega)}(u,x)+l_{(N,\omega)}(u,z)$. Moreover, $(N,\omega)$ is ultrametric, so we have $l_{(N,\omega)}(v,x)=l_{(N,\omega)}(v,y)$ and $l_{(N,\omega)}(u,x)=l_{(N,\omega)}(u,z)$. Finally, we have $l_{(N,\omega)}(u,x)=l_{(N,\omega)}(v,x)+\omega((u,v))$. In combination, these equalities imply $\omega((u,v))=\frac{1}{2}(\D(x,z)-\D(x,y))$. This concludes the proof that $\omega$ is uniquely determined by $N$ and $\D$, and thus the proof of the theorem.
\end{proof}

\section{Future directions}\label{sec:future}

In this paper we have introduced the concept of arboreal ultrametrics and
considered how to produce them from unrooted trees and how to characterize them. There remain several  interesting 
open directions for research on arboreal ultrametrics and related structures.

First, as stated in Corollary~\ref{cor:on4}, 
we can recognize whether or not a partial distance on a set of size $n$ is
an  arboreal ultrametric in $O(n^4)$ time. It would be interesting
to know if it is possible to obtain a better bound. 
Note that extending Bandelt's approach in \cite{bandelt1990recognition}
to recognizing whether or not a distance is an 
ultrametric  on a set of size $n$ in $O(n^2log(n))$ time
in the obvious way does appear to give a faster algorithm since 
checking Property~(U3) in Theorem~\ref{th:ultra} is problematic. Thus 
some other approach is probably required if it is indeed possible
to find an improvement (e.g. by considering approaches
such as those in \cite{culberson1989fast} which can recognize an ultrametric in $O(n^2)$ time).

More generally since distances can be defined in various way and for many different types of data \cite[Chapter 11]{Felsenstein2004},
it would be interesting to develop an algorithm that not 
only recognizes arboreal ultrametrics but also constructs them from distances
arising from real data.
The unweighted pair group method with arithmetic mean (UPGMA) \cite{sokal1958statistical} is
a popular algorithm that takes as input a distance matrix and outputs an
ultrametric tree. It also has the property that if the input is an ultrametric 
then is produces the unique ultrametric tree that realizes this ultrametric. It could
be worthwhile exploring if some algorithm could be developed for 
producing ultrametric arboreal networks from partial distances that generalizes UPGMA. 

In another but related direction, there are several results concerning
the approximation of distances by ultrametrics (see e.g. \cite[Section 3]{chepoi2000approximation}
for a review). It could be worth exploring which of these results might 
extend to approximations of partial distances by arboreal ultrametrics. For example,
it is well-known (see e.g. \cite{chepoi2000approximation}) that any distance has an ultrametric subdominant (or lower maximum approximation); is there such a result for partial distances? Note that
computing subdominants is closely related to the concept of the Farris tranform 
(see e.g. \cite{dress2007some} for a review), which could also be interesting
to investigate in the context of arboreal ultrametrics

Finally, as we have seen above, our characterization for
arboreal ultrametrics is closely related to the theory of symbolic ultrametrics.
In particular, for arboreal ultrametrics we are considering
the situation where the symbols are real numbers. 
It would be interesting to study what might happen if we replace 
real numbers with other algebraic structures such as groups. 
Note that this problem has already been considered 
in the context of symbolic ultrametrics (see e.g. \cite[Section 7.6]{SS03} for an overview). 
More generally, it was recently shown that 
the class of distance-hereditary graphs is precisely the class of 
undirected graphs that can be explained by arboreal networks \cite{S25}. Thus, 
it could be interesting to investigate if 
our results lead to new directions of study in 
the theory of distance-hereditary graphs.

\acknowledgments
The authors thank the reviewer for their comments.
GES thanks the University of East Anglia, for hosting him for two stays during which significant progress were made on this project.
Part of this work was done while the authors were in residence at the Institute for Computational and Experimental Research in Mathematics (ICERM) in Providence (RI, USA) during the \emph{Theory, Methods, and Applications of Quantitative Phylogenomics} program (supported by grant DMS-1929284 of the National Science Foundation (NSF)).

\bibliographystyle{abbrvnat}
\bibliography{bibliography}

\begin{thebibliography}{21}
\providecommand{\natexlab}[1]{#1}
\providecommand{\url}[1]{\texttt{#1}}
\expandafter\ifx\csname urlstyle\endcsname\relax
  \providecommand{\doi}[1]{doi: #1}\else
  \providecommand{\doi}{doi: \begingroup \urlstyle{rm}\Url}\fi

\bibitem[Bandelt(1990)]{bandelt1990recognition}
H.-J. Bandelt.
\newblock Recognition of tree metrics.
\newblock \emph{SIAM Journal on Discrete Mathematics}, 3\penalty0 (1):\penalty0
  1--6, 1990.

\bibitem[Buneman(1971)]{B71}
P.~Buneman.
\newblock The recovery of trees from measures of dissimilarity.
\newblock \emph{Mathematics in Archeological and Historical Sciences,
  Proceedings of the Anglo-Romanian Conference}, pages 387--395, 1971.

\bibitem[Chepoi and Fichet(2000)]{chepoi2000approximation}
V.~Chepoi and B.~Fichet.
\newblock $l_\infty$-approximation via subdominants.
\newblock \emph{Journal of Mathematical Psychology}, 44\penalty0 (4):\penalty0
  600--616, 2000.

\bibitem[Coleman et~al.(2015)Coleman, Hogan, McGoldrick, Rudden, and
  McInerney]{coleman2015}
O.~Coleman, R.~Hogan, N.~McGoldrick, N.~Rudden, and J.~McInerney.
\newblock Evolution by pervasive gene fusion in antibiotic resistance and
  antibiotic synthesizing genes.
\newblock \emph{Computation}, 3\penalty0 (2):\penalty0 114--127, 2015.

\bibitem[Culberson and Rudnicki(1989)]{culberson1989fast}
J.~C. Culberson and P.~Rudnicki.
\newblock A fast algorithm for constructing trees from distance matrices.
\newblock \emph{Information Processing Letters}, 30\penalty0 (4):\penalty0
  215--220, 1989.

\bibitem[D\"ocker et~al.(2019)D\"ocker, van Iersel, Kelk, and Linz]{D19}
J.~D\"ocker, L.~van Iersel, S.~Kelk, and S.~Linz.
\newblock Deciding the existence of a cherry-picking sequence is hard on two
  tree.
\newblock \emph{Discrete Applied Mathematics}, 260:\penalty0 131--143, 2019.

\bibitem[Dress et~al.(2007)Dress, Huber, and Moulton]{dress2007some}
A.~Dress, K.~T. Huber, and V.~Moulton.
\newblock Some uses of the farris transform in mathematics and
  phylogenetics—a review.
\newblock \emph{Annals of Combinatorics}, 11\penalty0 (1):\penalty0 1--37,
  2007.

\bibitem[Felsenstein(2004)]{Felsenstein2004}
J.~Felsenstein.
\newblock \emph{Inferring Phylogenies}.
\newblock Sinauer Associates, Sunderland, Massachusetts, 2004.
\newblock ISBN 9780878931774.

\bibitem[Haggerty et~al.(2014)Haggerty, Jachiet, Hanage, Fitzpatrick, Lopez,
  O’Connell, Pisani, Wilkinson, Bapteste, and McInerney]{HJHFLOPWBM13}
L.~Haggerty, P.-A. Jachiet, W.~Hanage, D.~Fitzpatrick, P.~Lopez,
  M.~O’Connell, D.~Pisani, M.~Wilkinson, E.~Bapteste, and J.~McInerney.
\newblock A pluralistic account of homology: adapting the models to the data.
\newblock \emph{Molecular biology and evolution}, 31\penalty0 (3):\penalty0
  501--516, 2014.

\bibitem[Howorka(1981)]{H81}
E.~Howorka.
\newblock A characterization of ptolemaic graphs.
\newblock \emph{Journal of Graph Theory}, 5\penalty0 (3):\penalty0 323--331,
  1981.

\bibitem[Huang and Knowles(2016)]{HK16}
H.~Huang and L.~Knowles.
\newblock Unforseen consequences of excluding missing data from next-generation
  sequences: simulation study of rad sequences.
\newblock \emph{Systematic Biology}, 65\penalty0 (3):\penalty0 357--365, 2016.

\bibitem[Huber et~al.(2022)Huber, Moulton, and Scholz]{HMS22}
K.~T. Huber, V.~Moulton, and G.~E. Scholz.
\newblock Forest-based networks.
\newblock \emph{Bulletin of Mathematical Biology}, 84\penalty0 (10):\penalty0
  119, 2022.

\bibitem[Huber et~al.(2024)Huber, Moulton, and Scholz]{huber2024shared}
K.~T. Huber, V.~Moulton, and G.~E. Scholz.
\newblock Shared ancestry graphs and symbolic arboreal maps.
\newblock \emph{SIAM Journal on Discrete Mathematics}, 38\penalty0
  (4):\penalty0 2553--2577, 2024.

\bibitem[Humphries et~al.(2013)Humphries, Linz, and Semple]{HLS13}
P.~J. Humphries, S.~Linz, and C.~Semple.
\newblock Cherry picking: A characterization of the temporal hybridization
  number for a set of phylogenies.
\newblock \emph{Bulletin of Mathematical Biology}, 75:\penalty0 1879--1890,
  2013.

\bibitem[Kettleborough et~al.(2015)Kettleborough, Dicks, Roberts, and
  Huber]{kettleborough2015reconstructing}
G.~Kettleborough, J.~Dicks, I.~N. Roberts, and K.~T. Huber.
\newblock Reconstructing (super) trees from data sets with missing distances:
  not all is lost.
\newblock \emph{Molecular Biology and Evolution}, 32\penalty0 (6):\penalty0
  1628--1642, 2015.

\bibitem[Kinene et~al.(2016)Kinene, Wainaina, Maina, and
  Boykin]{kinene2016rooting}
T.~Kinene, J.~Wainaina, S.~Maina, and L.~M. Boykin.
\newblock Rooting trees, methods for.
\newblock \emph{Encyclopedia of Evolutionary Biology}, page 489, 2016.

\bibitem[Scholz(2025)]{S25}
G.~E. Scholz.
\newblock Representing distance-hereditary graphs with multi-rooted trees.
\newblock \emph{Graphs and Combinatorics}, 41\penalty0 (111), 2025.

\bibitem[Scholz et~al.(2019)Scholz, Popescu, Taylor, Moulton, and
  Huber]{scholz2019osf}
G.~E. Scholz, A.-A. Popescu, M.~I. Taylor, V.~Moulton, and K.~T. Huber.
\newblock Osf-builder: a new tool for constructing and representing
  evolutionary histories involving introgression.
\newblock \emph{Systematic Biology}, 68\penalty0 (5):\penalty0 717--729, 2019.

\bibitem[Semple and Steel(2003)]{SS03}
C.~Semple and M.~Steel.
\newblock \emph{Phylogenetics}.
\newblock Oxford University Press on Demand, 2003.

\bibitem[Sokal and Michener(1958)]{sokal1958statistical}
R.~R. Sokal and C.~D. Michener.
\newblock A statistical method for evaluating systematic relationships.
\newblock \emph{University of Kansas Science Bulletin}, 38:\penalty0
  1409--1438, 1958.

\bibitem[Steel(2016)]{S16}
M.~Steel.
\newblock \emph{Phylogeny: Discrete and Random Processes in Evolution}.
\newblock SIAM, 2016.

\end{thebibliography}

\section*{Appendix: Proof of Proposition~\ref{pr:uur}}

\begin{proof}
Throughout the proof, we will use the notation given in Algorithm~\ref{alg:ulup}. Put $n=|X|$. Let $x_1, \ldots, x_n$ be the wcps computed at Line~\ref{l:wcps}. Let $(N_2,\omega_2)$ be the weighted network built at Lines~\ref{l:N0} and \ref{l:o0}, and for $i \in \{3, \ldots, n\}$, let $(N_i,\omega_i)$ be the weighted network obtained after the $(i-2)^{\text{th}}$ instance of the loop initiated at Line~\ref{l:fori}. Note that for $i \in \{2, \ldots, n\}$, $(N_i,\omega_i)$ is a weighted network on $\{x_1, \ldots, x_i\}$. Moreover, since the algorithm terminates immediately after the last instance of the loop, the weighted network $(N,\omega)$ returned by Algorithm~\ref{alg:ulup} satisfies $(N,\omega)=(N_n,\omega_n)$. Therefore, 
it suffices to show that for all $i \in \{2, \ldots, n\}$, $(N_i,\omega_i)$ is an ultrametric uprooting of $(T_i,\lambda_i)$.
We shall do this by using induction on $i$.

As base case, consider the case $i=2$. Then $T_2$ has two vertices $x_1$ and $x_2$, and one edge $\{x_1,x_2\}$. By construction, $N_2$ is the rooted tree with vertex set $\{r_2,x,y\}$ and arc set $\{(r_2,x_1),(r_2,x_2)\}$ (Line~\ref{l:N0}), and we have $\omega((r_2,x_1))=\omega((r_2,x_2))=\frac{1}{2}\lambda_2(\{x_1,x_2\})$ (Line~\ref{l:o0}). From there, one can easily verify that, $N_2$ is an uprooting of $T_2$, and since, $\omega((r_0,x_1))+\omega((r_0,x_2))=\lambda_2(\{x_1,x_2\})$, it follows that $(N_2,\omega_2)$ is a weight-preserving uprooting of $(T_2,\lambda_2)$. Moreover, $\omega((r_2,x_1))=\omega((r_2,x_2))$, so $(N_2,\omega_2)$ is an ultrametric arboreal network.

Now, let $i>2$, and suppose that $(N_{i-1},\omega_{i-1})$ is an ultrametric uprooting of $(T_{i-1},\lambda_{i-1})$. To ease notation, we now put $(N,\omega)=(N_{i-1},\omega_{i-1})$, and $(N^+,\omega^+)=(N_i,\omega_i)$ (that is, the network obtained after completion of the $(i-2)^{\text{th}}$ instance of the loop). We 
now show that $(N^+,\omega^+)$ is an ultrametric uprooting of $(T_i,\lambda_i)$, 
which will complete the proof of the proposition since it
implies that the network that is returned in Line~\ref{l:end} has the desired properties.

 Let $v_i$, $x_i$ and $x$ be the elements picked on Line~\ref{l:defvx}, and let $p_x$ be as defined at Line~\ref{l:defvp}. Note that 
 the operations inside of the loop starting on Line~\ref{l:fori} are  divided into two independent parts. The first part from Line~\ref{l:ifbig} to \ref{l:veq} is concerned with adding a vertex $v$ to the network $N$, and depends on the relative values of $\omega((p_x,x))$ and $\lambda_i(\{v_i,x\})$. The second part, from Line~\ref{l:ifdown} to Line~\ref{l:omtop}, is concerned with adding the leaf $x_i$ to the network $N$, and depends on the relative values of $\lambda_i(\{v_i,x_i\})$ and $l_{(N,\omega)}(v,z)$ for some $z \in C_N(v)$. Note that by our induction hypothesis, $(N,\omega)$ is ultrametric, so $l_{(N,\omega)}(v,z)$ is independent from the choice of $z$ in $C_N(v)$. In the following, we always understand $z$ as some element in $C_N(v)$, without mentioning this explicitly.

Next, note that the first part has three subcases, respectively $\omega((p_x,x))>\lambda_i(\{v_i,x\})$ (Line~\ref{l:ifbig}), $\omega((p_x,x))<\lambda_i(\{v_i,x\})$ (Line~\ref{l:ifsm}) and $\omega((p_x,x))=\lambda_i(\{v_i,x\})$ (Line~\ref{l:ifeq}), and the second part has two subcases, namely, $\lambda_i(\{v_i,x_i\})=l_{(N,\omega)}(v,z)$ (Line~\ref{l:ifdown}) and $\lambda_i(\{v_i,x_i\}) \neq l_{(N,\omega)}(v,z)$ (Line~\ref{l:iftop}). Therefore, we will complete the proof 
that $(N^+,\omega^+)$ is an ultrametric uprooting of $(T_i,\lambda_i)$
by considering each of the six resulting combinations.\\

\noindent{\bf (1):} $\omega((p_x,x))>\lambda_i(\{v_i,x\})$ and $\lambda_i(\{v_i,x_i\})=l_{(N,\omega)}(v,z)$ (Lines~\ref{l:ifbig} and \ref{l:ifdown}). In this case $N^+$ is obtained from $N$ by first subdividing $(p_x,x)$ with the introduction of a new vertex $v$ (Line~\ref{l:vbig}), and adding the arc $(v,x_i)$ (Line~\ref{l:vdown}). By our induction hypothesis, $N$ is an uprooting of $T_{i-1}$, so $\overline N$ and $T_{i-1}$ are isomorphic. Moreover, by construction, $\overline {N^+}$ is obtained from $\overline N$ by subdividing the edge incident to $x$ with the introduction of vertex $v$, and adding the edge $\{v,x_i\}$. Since $x$ and $x_i$ are part of a cherry in $T_i$, there are two possibilities. Either $v_i=u$ where $u$ is the vertex adjacent to $x$ in $T_{i-1}$, or $v_i$ has degree $3$ and is adjacent to $u$, $x$ and $x_i$.
The first possibility cannot hold since $\lambda_{i-1}(\{u,x\})=\lambda_i(\{v_i,x\})$
and $\omega((p_x,x)) \leq \lambda_{i-1}(\{u,x\})$ must hold, which contradicts our assumption that $\omega((p_x,x))>\lambda_i(\{v_i,x\})$. Therefore the 
second possibility holds, and so $v_i$ has degree $3$ and is adjacent to $u$, $x$ and $x_i$. In other words, $T_i$ can be obtained from $T_{i-1}$ by subdividing the edge $(u,x)$ with the introduction of vertex $v_i$, and adding the edge $\{v_i,x\}$. Since $T_{i-1}$ and $\overline N$ are isomorphic, it follows that $T_i$ and $\overline {N^+}$ are isomorphic. Hence, $N^+$ is an uprooting of $T_i$.

We now show that $(N^+,\omega^+)$ is an uprooting of $(T_i,\lambda_i)$.
Let $\{w,w'\}$ be an edge of $T_i$. Suppose first that neither $w=v_i$ nor $w'=v_i$ holds. In this case, $w$ and $w'$ are vertices of $N$, and by our induction hypothesis, we have $\lambda_{i-1}(\{w,w'\})=\overline \omega(\{w,w'\})$. Moreover, since $w$ and $w'$ are both distinct from $v_i$, we have  $\lambda_{i-1}(\{w,w'\})=\lambda_i(\{w,w'\})$ and $\overline {\omega^+}(\{w,w'\})=\overline \omega(\{w,w'\})$. Hence, $\overline {\omega^+}(\{w,w'\})=\lambda_i(\{w,w'\})$ follows. Suppose now that $w=v_i$. Then, $w'$ is one of $x,x_i$ or $u$. Since $\omega^+(\{v,x\})=\lambda_i(\{v_i,x\})$ (Line~\ref{l:ombig}), $\overline{\omega^+}(\{v_i,x\})=\lambda_i(\{v_i,x\})$ follows. Similarly, $\omega^+(\{v,x_i\})=\lambda_i(\{v_i,x_i\})$ (Line~\ref{l:omdown}), so $\overline{\omega^+}(\{v_i,x_i\})=\lambda_i(\{v_i,x_i\})$ follows. Finally, we have $\lambda_i(\{v_i,u\})=\lambda_{i-1}(\{u,x\})-\lambda_i(\{u,x\})$. If $p_x$ is not a root of $N$ of degree $2$, then $\overline{\omega^+}(\{v_i,u\})=\omega^+((p_x,v))=\omega((p_x,x))-\lambda_i(\{v_i,x\})$. By our induction hypothesis, $\omega((p_x,x))=\overline \omega \{u,x\}=\lambda_{i-1}(\{u,x\})$, and so $\overline{\omega^+}((v_i,u))=\lambda_i(\{v_i,u\})$. This concludes the proof that $(N^+,\omega^+)$ is an uprooting of $(T_i,\lambda_i)$.

It remains to show that $(N^+,\omega^+)$ is an ultrametric network. Let $w$ be a root of $N^+$, and let $z \in C_{N^+}(w)$. Note that $w \neq v$ holds, since $v$ is not a root of $N^+$. In particular, $w$ is also a root of $N$. If $z$ is distinct from $x$ and $x_i$, then $z$ is not a descendant of $v$. In particular, all arcs on the path from $w$ to $z$ in $N^+$ are arcs of $N$, and for all such arcs $a$, $\omega(a)=\omega^+(a)$ holds. Therefore, $l_{(N^+,\omega^+)}(w,z)=l_{(N,\omega)}(w,z)$. If $z=x$ then the fact that $w \neq v$, together with the fact that $\omega^+((p_x,v))+\omega^+((v,x))=\omega((p_x,x))$ (Line~\ref{l:ombig}), also implies that $l_{(N^+,\omega^+)}(w,z)=l_{(N,\omega)}(w,z)$. 
Using the induction hypothesis that $(N,\omega)$ is ultrametric, it follows that for all $z,z' \in C_N(w)$ distinct from $x_i$, $l_{(N^+,\omega^+)}(w,z)=l_{(N^+,\omega^+)}(w,z')$. Finally, consider the case where $z=x_i$. Then, we have $l_{(N^+,\omega^+)}(w,x_i)=l_{(N^+,\omega^+)}(w,v)+\omega^+((v,x_i))$. Since $\omega^+((v,x_i))=\omega((v,x))$ (Lines~\ref{l:ombig} and \ref{l:omdown}), it follows that $l_{(N^+,\omega^+)}(w,x_i)=l_{(N^+,\omega^+)}(w,x)$. Since, as shown above, $l_{(N^+,\omega^+)}(w,x)=l_{(N^+,\omega^+)}(w,z)$ holds for all $z \in C_{N^+}(w)$ distinct from $x_i$, we have $l_{(N^+,\omega^+)}(w,x_i)=l_{(N^+,\omega^+)}(w,z)$ for all $z \in C_{N^+}(w)$. This concludes the proof that $(N^+,\omega^+)$ is an ultrametric network, and therefore, the proof that $(N^+,\omega^+)$ is an ultrametric uprooting of $(T_i,\lambda_i)$ 
and the proof for Case {\bf (1)}.

\noindent{\bf (2):}
$\omega((p_x,x))>\lambda_i(\{v_i,x\})$ and $\lambda_i(\{v_i,x_i\}) \neq l_{(N,\omega)}(v,z)$ (Lines~\ref{l:ifbig} and \ref{l:iftop}). In this case $N^+$ is obtained from $N$ by first subdividing $(p_x,x)$ with the introduction of a new vertex $v$ (Line~\ref{l:vbig}), and adding a new vertex $r$ and the arcs $(r,x_i)$ and $(r,v)$ (Line \ref{l:vtop}).
By construction, $\overline {N^+}$ is obtained from $\overline N$ by subdividing the edge incident to $x$ with the introduction of vertex $v$, and adding the edge $\{v,x_i\}$. From there, arguments similar as in case (1) imply that $N^+$ is an uprooting of $T_i$.

To prove that $(N^+,\omega^+)$ is an uprooting of $(T_i,\lambda_i)$, we first show
that $\omega((r,v))>0$. By Line~\ref{l:omtop}, we have $\omega((r,v))=\frac{1}{2}(\lambda_i(\{v_i,x_i\})-l_{(N,\omega)}(v,z))$ for some $z \in C_N(v)$. Note that since by assumption, $\lambda_i(\{v_i,x_i\}) \neq l_{(N,\omega)}(v,z)$, we have $\omega((r,v)) \neq 0$. Therefore, it suffices to show that $\omega((r,v)) \geq 0$ holds. We next remark that the only possible choice for $z$ is $x$. So,we have $\omega((r,v))=\frac{1}{2}(\lambda_i(\{v_i,x_i\})-l_{(N,\omega)}(v,x))=\frac{1}{2}(\lambda_i(\{v_i,x_i\})-\omega((v,x)))=\frac{1}{2}(\lambda_i(\{v_i,x_i\})-\lambda_i(\{v_i,x\}))$, where the latter equality comes from Line~\ref{l:ombig}. Since $\{v_i,x_i\}$ is the long end of the cherry, $\lambda_i(\{v_i,x_i\}) \geq \lambda_i(\{v_i,x\})$ follows, from which we obtain $\omega((r,v)) \geq 0$. Hence, $\omega((r,v))>0$ holds as desired.

Now, let $\{w,w'\}$ be an edge of $T_i$. If neither $w=v_i$ nor $w'=v_i$ holds, then by the arguments of case (1), we have $\overline {\omega^+}(\{w,w'\})=\lambda_i(\{w,w'\})$. The same arguments imply that $\overline{\omega^+}(\{v_i,x\})=\lambda_i(\{v_i,x\})$ and $\overline{\omega^+}((v_i,u))=\lambda_i(\{v_i,u\})$ hold, as $\omega((p_x,v))$ and $\omega((v,x))$ are defined the same way as in case (1) (Line~\ref{l:ombig}). Finally we have $\overline{\omega^+}(\{v_i,x_i\})=\omega((r,v))+\omega((r,x_i))=\lambda_i(\{v_i,x_i\})$ (Line~\ref{l:omtop}). This concludes the proof that $(N^+,\omega^+)$ is an uprooting of $(T_i,\lambda_i)$.

It remains to show that $(N^+,\omega^+)$ is ultrametric. Let $w$ be a root of $N^+$. If $w=r$, then $C_{N^+}(w)$ contains exactly $2$ elements, $x$ and $x_i$. We have $l_{(N^+,\omega^+)}(r,x)=\omega((r,v))+\omega((v,x))=\frac{1}{2}(\lambda_i(\{v_i,x_i\})-l_{(N,\omega)}(v,z))+\lambda_i(\{v_i,x\})$ for some $z \in C_N(v)$ (Lines~\ref{l:ombig} and \ref{l:omtop}), and $l_{(N^+,\omega^+)}(r,x_i)=\omega((r,x_i))=\frac{1}{2}(\lambda_i(\{v_i,x_i\})+l_{(N,\omega)}(v,z))$ for some $z \in C_N(v)$ (Line~\ref{l:omtop}). Hence, $l_{(N^+,\omega^+)}(r,x)=l_{(N^+,\omega^+)}(r,x_i)$ holds. Suppose now that $w$ is distinct from $r$, and let $z \in C_{N^+}(w)$. 
Note that $z$ is distinct from $x_i$ since $w \neq r$. Now, let $a$ be an arc from $w$ to $z$. If $a$ is an arc of $N$, then $\omega(a)=\omega^+(a)$ holds. Otherwise, $a$ is one of the arcs $(p_x,v)$ or $(v,x)$. Next, we remark that since $v$ has indegree $2$ and outdegree $1$, and $w$ is distinct from $r$, the presence of one of these arcs on the path from $w$ to $z$ implies that the other is also present. Moreover, in that case, the arc $(p_x,x)$ is present on the path from $w$ to $z$ in $N$. Since $\omega((p_x,v))+\omega((v,x))=\omega((p_x,x))$ (Line~\ref{l:ombig}), $l_{(N^+,\omega^+)}(w,z)=l_{(N,\omega)}(w,z)$ always holds. Together with the induction hypothesis that $(N,\omega)$ is ultrametric, it follows that for all $z,z' \in C_N(w)$, $l_{(N^+,\omega^+)}(w,z)=l_{(N^+,\omega^+)}(w,z')$. This concludes the proof that $(N^+,\omega^+)$ is an ultrametric network, and so $(N^+,\omega^+)$ is an ultrametric uprooting of $(T_i,\lambda_i)$ which completes Case {\bf(2)}.

\noindent{\bf (3):}
$\omega((p_x,x))<\lambda_i(\{v_i,x\})$ and $\lambda_i(\{v_i,x_i\})=l_{(N,\omega)}(v,z)$ (Lines~\ref{l:ifsm} and \ref{l:ifdown}).  We claim that this case is impossible. Indeed, by definition, we have $l_{(N,\omega)}(v,z)=l_{(N,\omega)}(p_x,z)-\omega((p_x,v))$. Moreover, $\omega((p_x,v))=\lambda_i(\{v_i,x\})-\omega((p_x,x))>0$, where the equality comes from Line~\ref{l:omsm} and the inequality from our assumption $\omega((p_x,x))<\lambda_i(\{v_i,x\})$. Hence, we have $l_{(N,\omega)}(v,z)<l_{(N,\omega)}(p_x,z)$. Since $(N,\omega)$ is ultrametric by our induction hypothesis, we have $l_{(N,\omega)}(p_x,z)=l_{(N,\omega)}(p_x,x)$, and since $x$ is a child of $p_x$, $l_{(N,\omega)}(p_x,x)=\omega((p_x,x))$. By assumption, $\omega((p_x,x))<\lambda_i(\{v_i,x\})$. Putting these inequalities together, $l_{(N,\omega)}(v,z)<\lambda_i(\{v_i,x\})$ follows. Finally, since $\lambda_i(\{v_i,x_i\})=l_{(N,\omega)}(v,z)$ by assumption, we obtain $\lambda_i(\{v_i,x_i\})<\lambda_i(\{v_i,x\})$, which is impossible by our choice of $x_i$ as the long end of a cherry in $(T_i,\lambda_i)$.

\noindent{\bf (4):}
$\omega((p_x,x))<\lambda_i(\{v_i,x\})$ and $\lambda_i(\{v_i,x_i\}) \neq l_{(N,\omega)}(v,z)$ (Lines~\ref{l:ifsm} and \ref{l:iftop}). We first show that in this case, $p_x$ is a root of $N$ of outdegree $2$. So, suppose for contradiction that this is not the case. Then $p_x$ is a vertex of $T_{i-1}$, and therefore, also a vertex of $T_i$. In particular, we have $\lambda_{i-1}(\{p_x,x\}) \geq \lambda_i(\{v,x\})$. Since $\lambda_{i-1}(\{p_x,x\})=\omega((p_x,x))$, this contradicts our assumption that $\omega((p_x,x))<\lambda_i(\{v_i,x\})$.

We now show that $N^+$ is an uprooting of $T_i$.
Let $u$ be the child of $p_x$ distinct from $x$. Since $p_x$ is a root of $N$ of degree $2$ with children $u$ and $x$, and $T_i$ and $\overline N$ are isomorphic, we have that $u$ is the vertex adjacent to $x$ in $T_i$. 
We now consider two subcases, $u=v_i$ and $u \neq v_i$. First, if $u=v_i$, we have $\lambda_{i-1}(\{u,x\})=\lambda_i(\{v_i,x\})$, and since $\lambda_{i-1}(\{u,x\})=\overline \omega(\{u,x\})=\omega((p_x,u))+\omega((p_x,x))$, $\lambda_i(\{v_i,x\})=\omega((p_x,u))+\omega((p_x,x))$. Hence, we build $N^+$ from $N$ by putting $u=v$ (Line~\ref{l:vsm0}), adding a new vertex $r$, and adding the arcs $(r,v)$ and $(r,x_i)$ (Line~\ref{l:vtop}). In particular, $\overline {N^+}$ is obtained from $\overline N$ by adding the edge $\{u,x_i\}$. Since $\overline N$ and $T_{i-1}$ are isomorphic, and $T_i$ is obtained from $T_{i-1}$ by adding the arc $\{u,x_i\}$, it follows that $\overline {N^+}$ and $T_i$ are isomorphic. 
Second, if $u \neq v_i$, we have $\lambda_{i-1}(\{u,x\})>\lambda_i(\{v_i,x\})$, and since $\lambda_{i-1}(\{u,x\})=\overline \omega(\{u,x\})=\omega((p_x,u))+\omega((p_x,x))$, $\lambda_i(\{v_i,x\})<\omega((p_x,u))+\omega((p_x,x))$. Hence, we build $N^+$ from $N$ by subdividing the arc $(p_x,u)$ with the introduction of vertex $v$ (Line~\ref{l:vsm}), adding a new vertex $r$, and adding the arcs $(r,v)$ and $(r,x_i)$ (Line~\ref{l:vtop}). In particular, $\overline {N^+}$ is obtained from $\overline N$ by subdividing the edge $\{u,x\}$ with the introduction of vertex $v$, and adding the edge $\{v,x_i\}$. Since $\overline N$ and $T_{i-1}$ are isomorphic, and $T_i$ is obtained from $T_{i-1}$ by subdividin the arc $\{u,x\}$ with the introduction of vertex $v_i$, and adding the edge $\{v_i,x_i\}$, it follows that $\overline {N^+}$ and $T_i$ are isomorphic. In summary, in both cases, we can conclude that $N^+$ is an uprooting of $T_i$.

We now show that $(N^+,\omega^+)$ is an uprooting of $(T_i,\lambda_i)$.
First, we show that $\omega((r,v))>0$, and in case $\lambda_i(\{v_i,x\}) \neq \omega((p_x,u))+\omega((p_x,x))$, that $\omega((v,u))>0$. To see the latter, recall that in the last paragraph, we showed that $\lambda_i(\{v_i,x\}) \leq \omega((p_x,u))+\omega((p_x,x))$ always holds. In particular, $\omega((v,u))=\omega((p_x,u))+\omega((p_x,x))-\lambda_i(\{v_i,x\}) \geq 0$, and since $\lambda_i(\{v_i,x\}) \neq \omega((p_x,u))+\omega((p_x,x))$, the inequality is strict.
We now show that $\omega((r,v))>0$. By Line~\ref{l:omtop}, we have $\omega((r,v))=\frac{1}{2}(\lambda_i(\{v_i,x_i\})-l_{(N,\omega)}(v,z))$ for some $z \in C_N(v)$. By construction, and since $\omega((v,u))>0$, we have $l_{(N,\omega)}(v,z) \leq l_{(N,\omega)}(p_x,z)$. 
Moreover, $(N,\omega)$ is ultrametric by our induction hypothesis, so we have $l_{(N,\omega)}(p_x,z)=l_{(N,\omega)}(p_x,x)$. Since $l_{(N,\omega)}(p_x,x)=\omega((p_x,x))$, 
we have $l_{(N,\omega)}(v,z) \leq \omega((p_x,x))$. 
Moreover since by assumption, $\omega((p_x,x))<\lambda_i(\{v_i,x\})$, we get $l_{(N,\omega)}(v,z)<\lambda_i(\{v_i,x\})$. Finally, since $x_i$ is the long end of its cherry, we have $\lambda_i(\{v_i,x\}) \leq \lambda_i(\{v_i,x_i\})$, and so $l_{(N,\omega)}(v,z)<\lambda_i(\{v_i,x_i\})$ which proves, using the definition of $\omega((r,v))$ at Line~\ref{l:omtop}, that $\omega((r,v))>0$.

Now, let $\{w,w'\}$ be an edge of $T_i$. If $\{w,w'\}$ is an edge of $T_{i-1}$, then by the arguments 
used in Case {\bf (1)}, we have $\overline {\omega^+}(\{w,w'\})=\lambda_i(\{w,w'\})$. Suppose now that $\{w,w'\}$ is not an edge of $T_i$. If $v_i=u$, the only such edge is $\{v_i,x_i\}$. If otherwise, $v_i \neq u$, then there are three such edges, namely, $\{v_i,x_i\}$, $\{u,v_i\}$ and $\{v_i,x\}$. For $\{v_i,x_i\}$, we have $\overline{\omega^+}(\{v_i,x_i\})=\omega((r,v))+\omega((r,x_i))=\lambda_i(\{v_i,x_i\})$ (Line~\ref{l:omtop}). For $\{u,v_i\}$, we have $\overline{\omega^+}(\{u,v_i\})=\omega((v,u))=\omega((p_x,u))+\omega((p_x,x))-\lambda_i(\{v_i,x\})$ (Line~\ref{l:omsm}). Since $(N,\omega)$ is an uprooting of $(T_{i-1},\omega_{i-1})$, and $p_x$ is a root of $N$ of degree $2$, $\omega((p_x,u))+\omega((p_x,x))=\lambda_{i-1}(\{u,x\})$. Therefore, we have $\overline{\omega^+}(\{u,v_i\})=\lambda_{i-1}(\{u,x\})-\lambda_i(\{v_i,x\})=\lambda_i(\{u,v_i\})$. Finally, for $\{v_i,x\}$, we have $\overline{\omega^+}(\{v_i,x\})=\omega((p_x,x))+\omega((p_x,v))=\lambda_i(\{v_i,x\})$ (Line~\ref{l:omsm}). This concludes the proof that $(N^+,\omega^+)$ is an uprooting of $(T_i,\lambda_i)$.

It remains to show that $(N^+,\omega^+)$ is ultrametric. Let $w$ be a root of $N^+$. Suppose first that $w=r$, and let $z,z' \in C_N(w)$. If $z$ and $z'$ are distinct from $x_i$, then $z$ and $z'$ are both descendants of $v$. In particular, we have $l_{(N^+,\omega^+)}(r,z)=\omega((r,v))+ l_{(N^+,\omega^+)}(v,z)$ and $l_{(N^+,\omega^+)}(r,z')=\omega((r,v))+ l_{(N^+,\omega^+)}(v,z')$. By our induction hypothesis, $(N,\omega)$ is ultrametric, so we have $l_{(N^+,\omega^+)}(v,z)= l_{(N^+,\omega^+)}(v,z')$, and $l_{(N^+,\omega^+)}(r,z)= l_{(N^+,\omega^+)}(r,z')$. Otherwise, if $z'=x_i$, then we have $l_{(N^+,\omega^+)}(r,x_i) = \omega^+((r,x_i)) = \frac{1}{2}(\lambda_i(\{v_i,x_i\})+l_{(N,\omega)}(v,z))$ (Line~\ref{l:omtop}), and $l_{(N^+,\omega^+)}(r,z)= \omega((r,v))+l_{(N^+,\omega^+)}(v,z) =\frac{1}{2}(\lambda_i(\{v_i,x_i\})$ $- l_{(N,\omega)}(v,z)) + l_{(N^+,\omega^+)}(v,z)$. Hence $l_{(N^+,\omega^+)}(r,x_i)=l_{(N^+,\omega^+)}(r,z)$.  Suppose now that $w$ is distinct from $r$. Then, arguments similar to the ones used in Case {\bf (2)} (replacing arcs $(p_x,x)$ and $(v,x)$ with $(p_x,u)$ and $(v,u)$, respectively, and remarking that $\omega^+((p_x,v))+\omega^+((v,u))=\omega((p_x,u))$ holds by Line~\ref{l:omsm}) lead to the conclusion that $l_{(N^+,\omega^+)}(w,z)=l_{(N^+,\omega^+)}(w,z')$ holds for all $z,z' \in C_N(v)$.
It follows that $(N^+,\omega^+)$ is ultrametric, which completes Case~{\bf(4)}.

\noindent{\bf (5):}
$\omega((p_x,x))=\lambda_i(\{v_i,x\})$ and $\lambda_i(\{v_i,x_i\})=l_{(N,\omega)}(v,z)$ (Lines~\ref{l:ifeq} and \ref{l:ifdown}). This is a special case of {\bf (1)}, where $v$ is chosen as $p_x$ instead of resulting from a subdivision (Line~\ref{l:veq}).

\noindent{\bf (6):}
 $\omega((p_x,x))=\lambda_i(\{v_i,x\})$ and $\lambda_i(v_i,x_i) \neq l_{(N,\omega)}(v,z)$ (Lines~\ref{l:ifeq} and \ref{l:iftop}). This is 
 a special case of {\bf (2)}, where $v$ is chosen as $p_x$ instead of resulting from a subdivision (Line~\ref{l:veq}). \end{proof}


\end{document}